\documentclass[11pt]{article}
\usepackage{amsmath, upgreek}
\usepackage{amssymb}
\usepackage{amscd}
\usepackage{xspace}
\usepackage{verbatim}
\newcommand{\mysection}[1]{
\section{#1}\setcounter{equation}{0}}
\title{\bf Boundary singularities of solutions to elliptic viscous Hamilton-Jacobi equations}
\author{{\bf Tai Nguyen Phuoc}\quad
 {\bf Laurent V\'eron}\\[2mm]
{\small Laboratoire de Math\'ematiques et Physique Th\'eorique, }\\
{\small  Universit\'e Fran\c{c}ois Rabelais,  Tours,  FRANCE}}

\date{}
\begin{document}
 \maketitle


\newcommand{\txt}[1]{\;\text{ #1 }\;}
\newcommand{\tbf}{\textbf}
\newcommand{\tit}{\textit}
\newcommand{\tsc}{\textsc}
\newcommand{\trm}{\textrm}
\newcommand{\mbf}{\mathbf}
\newcommand{\mrm}{\mathrm}
\newcommand{\bsym}{\boldsymbol}
\newcommand{\scs}{\scriptstyle}
\newcommand{\sss}{\scriptscriptstyle}
\newcommand{\txts}{\textstyle}
\newcommand{\dsps}{\displaystyle}
\newcommand{\fnz}{\footnotesize}
\newcommand{\scz}{\scriptsize}
\newcommand{\be}{\begin{equation}}
\newcommand{\bel}[1]{\begin{equation}\label{#1}}
\newcommand{\ee}{\end{equation}}
\newcommand{\eqnl}[2]{\begin{equation}\label{#1}{#2}\end{equation}}
\newcommand{\barr}{\begin{eqnarray}}
\newcommand{\earr}{\end{eqnarray}}
\newcommand{\bars}{\begin{eqnarray*}}
\newcommand{\ears}{\end{eqnarray*}}
\newcommand{\nnu}{\nonumber \\}
\newtheorem{subn}{\name}
\renewcommand{\thesubn}{}
\newcommand{\bsn}[1]{\def\name{#1}\begin{subn}}
\newcommand{\esn}{\end{subn}}
\newtheorem{sub}{\name}[section]
\newcommand{\dn}[1]{\def\name{#1}}   
\newcommand{\bs}{\begin{sub}}
\newcommand{\es}{\end{sub}}
\newcommand{\bsl}[1]{\begin{sub}\label{#1}}
\newcommand{\bth}[1]{\def\name{Theorem}
\begin{sub}\label{t:#1}}
\newcommand{\blemma}[1]{\def\name{Lemma}
\begin{sub}\label{l:#1}}
\newcommand{\bcor}[1]{\def\name{Corollary}
\begin{sub}\label{c:#1}}
\newcommand{\bdef}[1]{\def\name{Definition}
\begin{sub}\label{d:#1}}
\newcommand{\bprop}[1]{\def\name{Proposition}
\begin{sub}\label{p:#1}}
\newcommand{\R}{\eqref}
\newcommand{\rth}[1]{Theorem~\ref{t:#1}}
\newcommand{\rlemma}[1]{Lemma~\ref{l:#1}}
\newcommand{\rcor}[1]{Corollary~\ref{c:#1}}
\newcommand{\rdef}[1]{Definition~\ref{d:#1}}
\newcommand{\rprop}[1]{Proposition~\ref{p:#1}}
\newcommand{\BA}{\begin{array}}
\newcommand{\EA}{\end{array}}
\newcommand{\BAN}{\renewcommand{\arraystretch}{1.2}
\setlength{\arraycolsep}{2pt}\begin{array}}
\newcommand{\BAV}[2]{\renewcommand{\arraystretch}{#1}
\setlength{\arraycolsep}{#2}\begin{array}}
\newcommand{\BSA}{\begin{subarray}}
\newcommand{\ESA}{\end{subarray}}
\newcommand{\BAL}{\begin{aligned}}
\newcommand{\EAL}{\end{aligned}}
\newcommand{\BALG}{\begin{alignat}}
\newcommand{\EALG}{\end{alignat}}
\newcommand{\BALGN}{\begin{alignat*}}
\newcommand{\EALGN}{\end{alignat*}}
\newcommand{\note}[1]{\textit{#1.}\hspace{2mm}}
\newcommand{\Proof}{\note{Proof}}
\newcommand{\qeda}{\hspace{10mm}\hfill $\square$}
\newcommand{\qed}{\\
${}$ \hfill $\square$}
\newcommand{\Remark}{\note{Remark}}
\newcommand{\modin}{$\,$\\[-4mm] \indent}
\newcommand{\forevery}{\quad \forall}
\newcommand{\set}[1]{\{#1\}}
\newcommand{\setdef}[2]{\{\,#1:\,#2\,\}}
\newcommand{\setm}[2]{\{\,#1\mid #2\,\}}
\newcommand{\mt}{\mapsto}
\newcommand{\lra}{\longrightarrow}
\newcommand{\lla}{\longleftarrow}
\newcommand{\llra}{\longleftrightarrow}
\newcommand{\Lra}{\Longrightarrow}
\newcommand{\Lla}{\Longleftarrow}
\newcommand{\Llra}{\Longleftrightarrow}
\newcommand{\warrow}{\rightharpoonup}
\newcommand{
\paran}[1]{\left (#1 \right )}
\newcommand{\sqbr}[1]{\left [#1 \right ]}
\newcommand{\curlybr}[1]{\left \{#1 \right \}}
\newcommand{\abs}[1]{\left |#1\right |}
\newcommand{\norm}[1]{\left \|#1\right \|}
\newcommand{
\paranb}[1]{\big (#1 \big )}
\newcommand{\lsqbrb}[1]{\big [#1 \big ]}
\newcommand{\lcurlybrb}[1]{\big \{#1 \big \}}
\newcommand{\absb}[1]{\big |#1\big |}
\newcommand{\normb}[1]{\big \|#1\big \|}
\newcommand{
\paranB}[1]{\Big (#1 \Big )}
\newcommand{\absB}[1]{\Big |#1\Big |}
\newcommand{\normB}[1]{\Big \|#1\Big \|}
\newcommand{\produal}[1]{\langle #1 \rangle}

\newcommand{\thkl}{\rule[-.5mm]{.3mm}{3mm}}
\newcommand{\thknorm}[1]{\thkl #1 \thkl\,}
\newcommand{\trinorm}[1]{|\!|\!| #1 |\!|\!|\,}
\newcommand{\bang}[1]{\langle #1 \rangle}
\def\angb<#1>{\langle #1 \rangle}
\newcommand{\vstrut}[1]{\rule{0mm}{#1}}
\newcommand{\rec}[1]{\frac{1}{#1}}
\newcommand{\opname}[1]{\mbox{\rm #1}\,}
\newcommand{\supp}{\opname{supp}}
\newcommand{\dist}{\opname{dist}}
\newcommand{\myfrac}[2]{{\displaystyle \frac{#1}{#2} }}
\newcommand{\myint}[2]{{\displaystyle \int_{#1}^{#2}}}
\newcommand{\mysum}[2]{{\displaystyle \sum_{#1}^{#2}}}
\newcommand {\dint}{{\displaystyle \myint\!\!\myint}}
\newcommand{\q}{\quad}
\newcommand{\qq}{\qquad}
\newcommand{\hsp}[1]{\hspace{#1mm}}
\newcommand{\vsp}[1]{\vspace{#1mm}}
\newcommand{\ity}{\infty}
\newcommand{\prt}{\partial}
\newcommand{\sms}{\setminus}
\newcommand{\ems}{\emptyset}
\newcommand{\ti}{\times}
\newcommand{\pr}{^\prime}
\newcommand{\ppr}{^{\prime\prime}}
\newcommand{\tl}{\tilde}
\newcommand{\sbs}{\subset}
\newcommand{\sbeq}{\subseteq}
\newcommand{\nind}{\noindent}
\newcommand{\ind}{\indent}
\newcommand{\ovl}{\overline}
\newcommand{\unl}{\underline}
\newcommand{\nin}{\not\in}
\newcommand{\pfrac}[2]{\genfrac{(}{)}{}{}{#1}{#2}}

\def\ga{\alpha}     \def\gb{\beta}       \def\gg{\gamma}
\def\gc{\chi}       \def\gd{\delta}      \def\ge{\epsilon}
\def\gth{\theta}                         \def\vge{\varepsilon}
\def\gf{\phi}       \def\vgf{\varphi}    \def\gh{\eta}
\def\gi{\iota}      \def\gk{\kappa}      \def\gl{\lambda}
\def\gm{\mu}        \def\gn{\nu}         \def\gp{\pi}
\def\vgp{\varpi}    \def\gr{\rho}        \def\vgr{\varrho}
\def\gs{\sigma}     \def\vgs{\varsigma}  \def\gt{\tau}
\def\gu{\upsilon}   \def\gv{\vartheta}   \def\gw{\omega}
\def\gx{\xi}        \def\gy{\psi}        \def\gz{\zeta}
\def\Gg{\Gamma}     \def\Gd{\Delta}      \def\Gf{\Phi}
\def\Gth{\Theta}
\def\Gl{\Lambda}    \def\Gs{\Sigma}      \def\Gp{\Pi}
\def\Gw{\Omega}     \def\Gx{\Xi}         \def\Gy{\Psi}

\def\CS{{\mathcal S}}   \def\CM{{\mathcal M}}   \def\CN{{\mathcal N}}
\def\CR{{\mathcal R}}   \def\CO{{\mathcal O}}   \def\CP{{\mathcal P}}
\def\CA{{\mathcal A}}   \def\CB{{\mathcal B}}   \def\CC{{\mathcal C}}
\def\CD{{\mathcal D}}   \def\CE{{\mathcal E}}   \def\CF{{\mathcal F}}
\def\CG{{\mathcal G}}   \def\CH{{\mathcal H}}   \def\CI{{\mathcal I}}
\def\CJ{{\mathcal J}}   \def\CK{{\mathcal K}}   \def\CL{{\mathcal L}}
\def\CT{{\mathcal T}}   \def\CU{{\mathcal U}}   \def\CV{{\mathcal V}}
\def\CZ{{\mathcal Z}}   \def\CX{{\mathcal X}}   \def\CY{{\mathcal Y}}
\def\CW{{\mathcal W}} \def\CQ{{\mathcal Q}}
\def\BBA {\mathbb A}   \def\BBb {\mathbb B}    \def\BBC {\mathbb C}
\def\BBD {\mathbb D}   \def\BBE {\mathbb E}    \def\BBF {\mathbb F}
\def\BBG {\mathbb G}   \def\BBH {\mathbb H}    \def\BBI {\mathbb I}
\def\BBJ {\mathbb J}   \def\BBK {\mathbb K}    \def\BBL {\mathbb L}
\def\BBM {\mathbb M}   \def\BBN {\mathbb N}    \def\BBO {\mathbb O}
\def\BBP {\mathbb P}   \def\BBR {\mathbb R}    \def\BBS {\mathbb S}
\def\BBT {\mathbb T}   \def\BBU {\mathbb U}    \def\BBV {\mathbb V}
\def\BBW {\mathbb W}   \def\BBX {\mathbb X}    \def\BBY {\mathbb Y}
\def\BBZ {\mathbb Z}

\def\GTA {\mathfrak A}   \def\GTB {\mathfrak B}    \def\GTC {\mathfrak C}
\def\GTD {\mathfrak D}   \def\GTE {\mathfrak E}    \def\GTF {\mathfrak F}
\def\GTG {\mathfrak G}   \def\GTH {\mathfrak H}    \def\GTI {\mathfrak I}
\def\GTJ {\mathfrak J}   \def\GTK {\mathfrak K}    \def\GTL {\mathfrak L}
\def\GTM {\mathfrak M}   \def\GTN {\mathfrak N}    \def\GTO {\mathfrak O}
\def\GTP {\mathfrak P}   \def\GTR {\mathfrak R}    \def\GTS {\mathfrak S}
\def\GTT {\mathfrak T}   \def\GTU {\mathfrak U}    \def\GTV {\mathfrak V}
\def\GTW {\mathfrak W}   \def\GTX {\mathfrak X}    \def\GTY {\mathfrak Y}
\def\GTZ {\mathfrak Z}   \def\GTQ {\mathfrak Q}

\font\Sym= msam10 
\def\SYM#1{\hbox{\Sym #1}}
\newcommand{\bdw}{\prt\Gw\xspace}
\tableofcontents
\date{}
\maketitle\medskip

\noindent{\small {\bf Abstract} We study the boundary value problem with measures for (E1) $-\Gd u+g(|\nabla u|)=0$ in a bounded domain $\Gw$ in $\BBR^N$, satisfying (E2) $ u=\gm$ on $\prt\Gw$ and prove that if $g\in L^1(1,\infty;t^{-(2N+1)/N}dt)$  is nondecreasing (E1)-(E2) can be solved with any positive bounded measure. When $g(r)\geq r^q$ with $q>1$ we prove that any positive function satisfying (E1) admits a boundary trace which is an outer regular Borel measure, not necessarily bounded. When $g(r)=r^q$ with $1<q<q_c=\frac{N+1}{N}$ we prove the existence of a positive solution with a general outer regular Borel measure $\gn\equiv\!\!\!\!\!\!/\;\infty$ as boundary trace and characterize the boundary isolated singularities of positive solutions. When  $g(r)=r^q$ with $q_c\leq q<2$ we prove that a necessary condition for solvability is that $\gm$ must be absolutely continuous with respect to the Bessel capacity $C_{\frac{2-q}{q},q'}$. We also characterize boundary removable sets for moderate and sigma-moderate solutions. 

\vspace{1mm}
\hspace{.05in}
\medskip
\mysection{Introduction}
Let  $\Gw \sbs \BBR^N$ be a bounded domain with $C^2$ boundary and $g:\BBR_+ \to\BBR_+$ a nondecreasing continuous function vanishing at $0$. In this article we investigate several boundary data questions associated to nonnegative solutions of the following equation
\bel{B*} -\Gd u + g(\abs{\nabla u}) = 0 \qq \text{in } \Gw,
\ee
and we emphasize on the particular case of 
\bel{B} -\Gd u + \abs{\nabla u}^q = 0 \qq \text{in } \Gw.
\ee
where $q$ is a real number mainly in the range $1<q<2$. We investigate first the generalized boundary value problem with measure associated to $(\ref{B*})$
\bel{B1}\left\{\BA {l} -\Gd u + g(\abs{\nabla u}) = 0 \qq \text{in } \Gw\\[1mm]
\phantom{-\Gd  +g(\abs{\nabla u})}
u=\gm \qq \text{on } \prt\Gw
\EA\right.\ee
where $\gm$ is a measure on $\prt\Gw$. By a solution we mean an integrable function $u$ such that $g(|\nabla u|)\in L^1_d(\Gw)$ where $d=d(x):=\dist (x,\prt\Gw)$ satisfying
\bel{B2}
\myint{\Gw}{}\left(-u\Gd\gz+g(\abs{\nabla u})\gz\right)dx=-\myint{\prt\Gw}{}\myfrac{\prt \gz}{\prt \bf n}d\gm
\ee
for all  $\gz\in X(\Gw):=\{\gf\in C^{1}_0(\overline\Gw):\Gd\gf\in L^{\infty}(\Gw)\}$, where $\bf n$ denotes the normal outward unit vector to $\prt \Gw$.  The {\it  integral subcriticality condition} for $g$ is the following
\bel{B2'}
\myint{1}{\infty}g(s)s^{-\frac{2N+1}{N}}ds<\infty.
\ee
When $g(r)\leq r^q$, this condition is satisfied if $0<q<q_c:=\frac{N+1}{N}$. Our main existence result is the following.\smallskip

\bth{Bth1} Assume $g$ satisfies $(\ref{B2'})$. Then for any positive bounded Borel measure $\gm$ on $\prt\Gw$ there exists a maximal positive solution $\ovl u_\gm$ to problem $(\ref{B1})$. Furthermore the problem is closed for weak convergence of boundary data.
\es

Note that we do not know if problem $(\ref{B2})$ has a unique solution, {\it except if $g(r)=r^q$ with $0<q<q_c$ and $\gm=c\gd_0$ in which case we prove that uniqueness holds}. 
A natural way for studying $(\ref{B*})$ is to introduce the notion of {\it boundary trace}. When $g(r)\geq r^q$ with $q>1$ we prove in particular that the following result holds in which statement we denote $\Gs_\gd=\{x\in\Gw:d(x)=\gd\}$ for $\gd>0$:

\bth{Bth2} Let $u$ be any positive solution of $(\ref{B*})$. Then for any $x_0\in\prt\Gw$ the following dichotomy occurs:\smallskip

\noindent(i) Either there exists an open neighborhood $U$ of $x_0$ such that
\bel{BB0} \myint{\Gw\cap U}{}g(|\nabla u|)d(x)dx<\infty
\ee
and there exists a positive Radon measure $\gm_U$ on $\prt\Gw\cap U$ such that $u|_{_{\Gs_\gd \cap U}}$ converges to $\gm_U$ in the weak sense of measures when $\gd\to 0$.\smallskip

\noindent(ii) Or for any open neighborhood $U$ of $x_0$ there holds
\bel{BB1} \myint{\Gw\cap U}{}g(|\nabla u|)d(x)dx=\infty,
\ee
and
\bel{BB2} \lim_{\gd\to 0}\myint{\Gs_\gd\cap U}{}udS=\infty.
\ee
\es

The set $\CS(u)$ of boundary points $x_0$ with the property (ii) is closed and there exists a unique Borel measure $\gm$ on $\CR(u):=\prt\Gw\setminus \CS(u)$ such that $u|_{_{\Gs_\gd}}$ converges to $\gm$ in the weak sense of measures on $\CR(u)$. The couple $(\CS(u),\gm)$ is the boundary trace of $u$, denoted by $tr_{\prt\Gw}(u)$. The trace  framework has also the advantage of pointing out some of the main questions which remain to be solved as it was done for the semilinear equation
\bel{A0} -\Gd u + h(u) = 0 \qq \text{in } \Gw.
\ee
and the associated Dirichlet problem with measure
\bel{A1}\left\{\BA {l} -\Gd u + h(u) = 0 \qq \text{in } \Gw\\[1mm]
\phantom{-\Gd  + h(u)}
u=\gm \qq \text{on } \prt\Gw,
\EA\right.\ee
where $h:\BBR\to\BBR$ is a continuous nondecreasing function vanishing at $0$. Much is known since the first paper of Gmira and V\'eron \cite {GV} and many developments are due to Marcus and V\'eron \cite{MV1}--\cite{MV4}  in particular when 
 $(\ref{A0})$ is replaced by 
\bel{A2} -\Gd u + \abs{u}^{q-1}u = 0 \qq \text{in } \Gw.
\ee
with $q>1$. We recall below some of the main aspects of the results dealing with $(\ref{A0})$--$(\ref{A2})$, this will play the role of the breadcrumbs trail for our study. \smallskip

- Problem $(\ref{A1})$ can be solved (in a unique way) for any bounded measure $\gm$ if $h$ satisfies

\bel{A1+1}\myint{1}{\infty}(h(s)+|h(-s)|)s^{-\frac{2N}{N-1}}ds<\infty.
\ee 
If $h(u)= \abs{u}^{q-1}u$ the condition $(\ref{A1+1})$ is verified if and only if $1<q<q_s$, {\it the subcritical range}; $q_s=\frac{N+1}{N-1}$ is a critical exponent for $(\ref{A2})$.

- When $1<q<q_s$, boundary isolated singularities of nonnegative solutions of $(\ref{A2})$ can be completely characterized i.e. if  $u\in C(\overline\Gw\setminus\{0\})$ is a nonnegative solution of $(\ref{A2})$ vanishing on $\prt\Gw\setminus\{0\}$, then either it solves the associated Dirichlet problem with $\gm=c\gd_0$ for some $c\geq 0$ ({\it weak singularity}), or
\bel{A4}
u(x)\approx d(x)|x|^{-\frac{q+1}{q-1}}\qq\text{as }x\to 0. \q \text{({\it strong singularity})}
\ee

- Always in the subcritical range  it is proved that for any couple $(\CS,\gm)$ where $\CS\subset\prt\Gw$ is closed and $\gm$ is a positive Radon measure on $\CR=\prt\Gw\setminus \CS$ there exists a unique positive solution $u$ of $(\ref{A2})$ with boundary trace $(\CS,\gm)$ (in the sense defined in \rth {Bth2}).

- When $q\geq q_s$, i.e. the {\it supercritical range}, any solution $u\in C(\overline\Gw\setminus\{0\})$ of $(\ref{A2})$ vanishing on $\prt\Gw\setminus\{0\}$ is identically $0$, i.e. {\it isolated boundary singularities are removable}. This result due to Gmira-V\'eron has been extended, either by probabilistic tools by Le Gall \cite {Lg1}, \cite {Lg2}, Dynkin \cite{Dy1}, Dynkin and Kuznetsov \cite {DK1}, \cite {DK2}, with the restriction $ q_s\leq q\leq 2$, or by purely analytic methods by Marcus and V\'eron \cite {MV1}, \cite {MV2} in the whole range $ q_s\leq q$. The key tool for describing the problem is the Bessel capacity $C_{\frac{2}{q},q'}$ in dimension $N-1$ (see \cite{AdHe} for a detailled presentation of capacities). We list some of the most striking results. The associated Dirichlet problem can be solved with $\gm\in\mathfrak M^+(\prt\Gw)$ if and only if $\gm$ is absolutely continuous with respect to the $C_{\frac{2}{q},q'}$-capacity. If $K\subset\prt\Gw$ is compact and $u\in C(\overline\Gw\setminus K)$ is a solution of $(\ref{A2})$ vanishing on $\prt \Gw\setminus K$, then $u$ is necessary zero if and only if $C_{\frac{2}{q},q'}(K)=0$. The complete characterization of positive solutions of $(\ref{A2})$ has been obtained by Mselati \cite {Ms} when $q=2$, Dynkin \cite {Dy2} when $ q_s\leq q\leq 2$, and finally Marcus \cite {Ma} when $ q_s\leq q$; they proved in particular that any positive solution $u$ is {\it sigma-moderate}, i.e. that there exists an increasing sequence of positive measures $\gm_n \in\mathfrak M^+(\prt\Gw)$ such that the sequence of the solutions $u=u_{\gm_n}$ of the associated Dirichlet  problem with $\gm=\gm_n$ converges to $u$. \medskip

Concerning $(\ref{B})$ we prove an existence result of solutions with a given trace belonging to the class of general outer regular Borel measures (not necessarily locally bounded).
\bth{Bth1+1} Assume $1<q<q_c$ and $\CS\subsetneq\prt\Gw$ is closed and $\gm$ is a positive Radon measure on $\CR:=\prt\Gw\setminus\CS$, then there exists a positive solution $u$ of $(\ref{B})$ such that $tr_{\prt\Gw}(u)=(\CS,\gm)$. 
\es

When $1<q<q_c$ we prove a stronger result, using the characterization of singular solutions with strong singularities (see \rth{Bth3+1} below). When $q_c\leq q<2$ we prove that \rth{Bth1+1} still holds with $\gm=0$ if $\CS=\overline G$ where $G\subsetneq\prt\Gw$ is relatively open, $\prt G$ satisfies an interior sphere condition. Surprisingly the condition $\CS\subsetneq\prt\Gw$ is necessary since there cannot exists any {\it large solution}, i.e. a solution which blows-up everywhere on $\prt\Gw$. 
 \medskip

In order to characterize isolated singularities of positive solutions of $(\ref{B})$ we introduce the following problem on the upper hemisphere $S^{N-1}_+$ of the unit sphere in $\BBR^N$
\bel{B3}\left\{\BA {l}
-\Gd'\gw+\left(\left(\frac{2-q}{q-1}\right)^2\gw^2+|\nabla'\gw|^2\right)^{\frac{q}{2}}-\frac{2-q}{q-1}\left(\frac{q}{q-1}-N\right)\gw
=0\q\text{in }S^{N-1}_+\\[2mm]
\phantom{-\Gd'+\left(\left(\frac{2-q}{q-1}\right)^2\gw^2+|\nabla'\gw|^2\right)^{\frac{q}{2}}-\frac{2-q}{q-1}\left(\frac{q}{q-1}-N\right)\gw}
\gw=0\q\text{on }\prt S^{N-1}_+,
\EA\right.\ee
where $\nabla'$ and $\Gd'$ denote respectively the covariant gradient and the Laplace-Beltrami operator on $S^{N-1}$. To  any solution $\gw$ of $(\ref{B3})$ we can associate a singular separable solution $u_s$ of $(\ref{B})$ in $\BBR^N_+:=\{x=(x_1,x_2,...,x_N)=(x',x_N):x_N>0\}$ vanishing on $\prt \BBR^N_+ \sms \{0\}$ written in spherical coordinates $(r,\gs)=(|x|,\frac{x}{|x|})$
\bel{B4}
u_s(x)=u_s(r,\gs)=r^{-\frac{2-q}{q-1}}\gw(\gs)\qq\forall x\in \overline{\BBR^N_+}\setminus\{0\}.
\ee

\bth{Bth2+1} Problem $(\ref{B3})$ admits a positive solution if and only if $1<q<q_c$. Furthermore this solution is unique and denoted by $\gw_s$. 
\es

This singular solution plays a fundamental role for describing isolated singularities.
\bth{Bth3} 
Assume $1<q<q_c$ and $u\in C^2(\Gw)\cap C(\overline\Gw\setminus\{0\})$ is a nonnegative solution of $(\ref{B})$ which vanishes on $\prt\Gw\setminus\{0\}$. Then the following dichotomy occurs:\smallskip

\noindent (i) Either there exists $c\geq 0$ such that $u=u_{c\gd_0}$ solves $(\ref{B1})$ with $g(r)=r^q$, $\gm=c\gd_0$ and 
\bel{B5}
u(x)=cP^\Gw(x,0)(1+o(1))\qq\text{as }x\to 0
\ee
where $P^\Gw$ is the Poisson kernel in $\Gw$.
\smallskip

\noindent (ii) Or $u=\lim_{c\to\infty}u_{c\gd_0}$ and 
\bel{B6}
\lim_{\tiny\BA{c}\Gw \ni x\to 0\\
\frac{x}{|x|}=\gs\in S^{N-1}_+
\EA}|x|^{\frac{2-q}{q-1}}u(x)=\gw_s(\gs).
\ee
\es

We also give a sharp estimate from below for singular points of the trace
\bth{Bth3+1} 
Assume $1<q<q_c$ and $u$ is a positive solution of $(\ref{B})$ with boundary trace $(\CS(u),\gm)$. Then for any $z\in\CS(u)$ there holds
\bel{B5+1}
u(x)\geq u_{\infty\gd_z}(x):=\lim_{c\to\infty}u_{c\gd_z}(x)\qq\forall x\in\Gw.
\ee
The description of $u_{\infty\gd_z}$ is provided by $u_s$ defined in $(\ref{B4})$, up to a translation and a rotation.
\es

The critical exponent $q_c$ plays for $(\ref{B})$ a role  similar to that of $q_s$ plays for $(\ref{A2})$ which is a consequence of the following theorem

\bth{Bth4} Assume $q_c \leq q <2$, then any nonnegative solution $u\in C^2(\Gw)\cap C(\overline\Gw\setminus\{0\})$ of $(\ref{B})$ vanishing on $\prt\Gw\setminus\{0\}$ is identically zero. 
\es

The supercritical case for equation $(\ref{B})$ can be understood using the Bessel capacity $C_{\frac{2-q}{q},q'}$ in dimension $N-1$, however we can only deal with moderate and sigma-moderate solutions. Following Dynkin \cite{Dy2}, \cite{DK3} we define

\bdef {mod}A positive solution $u$ of $(\ref{B})$ is moderate if there exists a bounded Borel measure $\gm\in \mathfrak M^+(\prt \Gw)$ such that $u$ solves problem $(\ref{B1})$ with $g(r)=r^q$. It is sigma-moderate if there exists an increasing sequence of solutions $\{u_{\gm_n}\}$, with boundary data $\{\gm_n\}\in \mathfrak M^+(\prt \Gw)$, which converges to $u$ when $n\to\infty$, locally uniformly in $\Gw$.
\es

Notice that the boundary trace theorem implies that the sequence $\{\gm_n\}$ is increasing. Equivalently we shall prove that a positive solution $u$ is moderate if and only if it is integrable in $\Gw$ and $\abs{\nabla u} \in L^q_d(\Gw)$. 

\bth{Bth5} Assume $q_c \leq q<2$ and $K\subset\prt\Gw$ is compact and satisfies $C_{\frac{2-q}{q},q'}(K)=0$. Then any positive moderate solution $u\in C^2(\Gw)\cap C(\overline\Gw\setminus K)$ of $(\ref{B})$ vanishing on $\prt\Gw\setminus K$ is identically zero. 
\es

As a corollary we prove that the above result remains true if $u$ is a sigma-moderate solution of $(\ref{B})$. The counterpart of this result is the following necessary condition for solving problem $(\ref{B1})$.

\bth{Bth6} Assume $q_c \leq q <2$ and $u$ is a positive moderate solution of $(\ref{B})$ with boundary data $\gm\in \mathfrak M^+(\prt\Gw)$. Then $\gm$ is absolutely continuous with respect to the $C_{\frac{2-q}{q},q'}$-capacity.
\es

For the sake of completeness we give, in Section 5, the results corresponding to the two extreme cases, $q=2$ and $q=1$ for equation $(\ref{B})$. If $q=2$ the Hopf-Cole change of unknown $u=\ln v$ transforms $(\ref{B})$ into a Poisson equation. When  $q=1$, equation $(\ref{B})$ is homogeneous of order $1$ and the equation inherits many properties of the Laplace equation.\medskip 

We end this article with a result concerning the question of existence and removability of solutions of
\bel{IA1}
-\Gd u+g(|\nabla u|)=\gm\qq\text{in }\Gw
\ee
where $\Gw$ is a bounded domain in $\BBR^N$ and $\gm$ a positive bounded Radon measure on $\Gw$. We prove that if $g$ is a locally Lipschitz nondecreasing function vanishing at $0$ and such that 
\bel{IA2}
\myint{1}{\infty}g(s)s^{-\frac{2N-1}{N-1}}ds<\infty
\ee
then problem $(\ref{IA1})$ admits a solution. In the power case

\bel{IA3}
-\Gd u+|\nabla u|^q=\gm\qq\text{in }\Gw
\ee
with $1<q<2$, the critical exponent is $q^*=\frac{N}{N-1}$. We prove that a necessary condition for solving $(\ref{IA3})$ with a positive Radon measure $\gm$  is that $\gm$ vanishes on Borel subsets $E$ with $C_{1,q'}$-capacity zero. The associated removability statement asserts that if $K$ a compact subset of $\Gw$ such that $C_{1,q'}(K)=0$, any positive solution of 
\bel{IA4}
-\Gd u+|\nabla u|^q=0\qq\text{in }\Gw\setminus K
\ee
 is bounded and can be extended as a solution to the whole $\Gw$.
\mysection{The Dirichlet problem and the boundary trace}

Throughout this article $\Gw$ is a bounded domain in $\BBR^N$ ($N\geq 2$) with a $C^2$ boundary $\prt\Gw$ and $c$ will denote a positive constant, independent of the data, the value of which may change from line to line. When needed the constant will be denoted by  $c_i$ or $C_i$ for some indices $i=1,2,...$, or some dependence will be made explicit such as $c(a,b,...)$ for some data $a$, $b$...For $r>0$ and $x \in \BBR^N$, we denote by $B_r(x)$ the ball with radius $r$ and center $x$. If $x=0$ we write $B_r$ instead of $B_r(0)$.
\subsection{Boundary data bounded measures}
We consider the following problem where $\gm$ belongs to the set $\GTM(\prt\Gw)$ of bounded Borel measures on $\prt\Gw$
\bel{M1}\left\{\BA{ll} -\Gd u+g(|\nabla u|)=0\qq\text{in }\Gw\\[2mm]\phantom{-\Gd +g(|\nabla u|)}
u=\gm\qq\text{on }\prt\Gw.\EA\right.
\ee
We assume that $g$ belongs to the class $\CG_0$ which means that $g:\BBR_+\to \BBR_+$ is a locally Lipschitz continuous nonnegative and nondecreasing function vanishing at $0$. The {\it integral subcriticality condition}  is the following
\bel{M}
\myint{1}{\infty}g(s)s^{-\frac{2N+1}{N}}ds<\infty.
\ee
 If $g(r)=r^q$ the integral subcriticality condition is satisfied if $0<q<q_c:=\frac{N+1}{N}$.

\bdef{sol} A function $u\in L^1(\Gw)$ such that $g(|\nabla u|)\in L^1_d(\Gw)$ is a weak solution of $(\ref{M1})$ if
\bel{M2}
\myint{\Gw}{}\left(-u\Gd\gz+g(|\nabla u|)\gz\right)dx=-\myint{\prt\Gw}{}\myfrac{\prt \gz}{\prt \bf n}d\gm
\ee
for all $\gz\in X(\Gw):=\{\gf\in C^1_0(\overline\Gw):\Gd\gf\in L^\infty(\Gw)\}$.
\es

If we denote respectively by $G^\Gw$ and $P^\Gw$  the Green kernel and the Poisson kernel in $\Gw$, with corresponding operators $\BBG^\Gw$ and $\BBP^\Gw$ it is classical from linear theory that the above definition is equivalent to

\bel{M2'}
u=\BBP^\Gw[\gm]-\BBG^\Gw[g(|\nabla u|)].
\ee

We recall that $M^p_h(\Gw)$ denote the Marcinkiewicz space (or weak $L^p$ space) of exponent $p\geq 1$ and weight $h>0$ defined by
\bel{Ma}
M^p_h(\Gw)=\left\{v\in L^1_{loc}(\Gw): \exists C\geq 0 \text{ s. t. } \myint{E}{}|v|hdx\leq C|E|_h^{1-\frac{1}{p}},\forall E\subset \Gw, E\text{ Borel}\right\},
\ee
where $|E|_h=\int_{}\chi_{_{E}} hdx$. The smallest constant $C$ for which $(\ref{Ma})$ holds is the Marcinkiewicz quasi-norm of $v$ denoted by
$\norm{v}_{M^p_h(\Gw)}$ and the following inequality will be much useful:
\bel{Ma'}
|\{x:|v(x)|\geq \gl\}|_h\leq \gl^{-p}\norm{v}_{M^p_h(\Gw)}^p\qq\forall \gl>0.
\ee

The main result of this section is the following existence and stability result  for problem $(\ref{M1})$.

\bth{meas} Assume $g\in\CG_0$ satisfies $(\ref{M})$, then for any $\gm\in\GTM^+(\prt\Gw)$ there exists a maximal solution $\bar u=\bar u_\gm$ to problem $(\ref {M1})$. Furthermore $\bar u\in M^{\frac{N}{N-1}}(\Gw)$ and $\abs{\nabla \bar u} \in M^{\frac{N+1}{N}}_d(\Gw)$. Finally, if $\{\gm_n\}$ is a sequence of positive bounded measures on $\prt\Gw$ which converges to $\gm$ in the weak sense of measures and $\{u_{\gm_n}\}$ is a sequence of solutions of $(\ref {M1})$ with boundary data $\gm_n$, then there exists a subsequence such that $\{ u_{\gm_{n_k}}\}$ converges to a solution $ u_\gm$ of $(\ref {M1})$ in $L^1(\Gw)$ and $\{g(|\nabla  u_{\gm_{n_k}}|)\}$ converges to $g(|\nabla u_{\gm}|)$ in $L^1_d(\Gw)$.
\es

We recall the following estimates \cite{BVi}, \cite{GV}, \cite{V1} and \cite{V2}.
\bprop{E1} For any $\ga\in [0,1]$, there exist a positive constant $c_1$ depending on $\ga$,  $\Gw$ and $N$ such that 
\bel{E1-1}
\norm{\BBG^{\Gw}[\gn]}_{L^1(\Gw)}+\norm{\BBG^{\Gw}[\gn]}_{M_{d^\ga}^{\frac{N+\ga}{N+\ga-2}}(\Gw)}\leq c_1\norm \gn_{\mathfrak M_{d^\ga}(\Gw)},
\ee
\bel{E1-2}
\norm{\nabla\BBG^{\Gw}[\gn]}_{M_{d^\ga}^{\frac{N+\ga}{N+\ga-1}}(\Gw)}\leq c_1\norm \gn_{\mathfrak M_{d^\ga}(\Gw)},
\ee
where
\bel{new}
\norm \gn_{\mathfrak M_{d^\ga}(\Gw)}:=\myint{\Gw}{}d^\ga(x) d|\gn| \qq\forall \gn\in \mathfrak M_{d^\ga}(\Gw),
\ee
\bel{E1-3}
\norm{\BBP^{\Gw}[\gm]}_{L^1(\Gw)}+\norm{\BBP^{\Gw}[\gm]}_{M^{\frac{N}{N-1}}(\Gw)}+\norm{\BBP^{\Gw}[\gm]}_{M_d^{\frac{N+1}{N-1}}(\Gw)}\leq c_1\norm \gm_{\GTM(\prt \Gw)},
\ee
\bel{E1-4}
\norm{\nabla\BBP^{\Gw}[\gm]}_{M_d^{\frac{N+1}{N}}(\Gw)}\leq c_1\norm \gm_{\GTM(\prt \Gw)},
\ee
for any $\gn\in \GTM_{d^\ga}(\Gw)$ and any $\gm\in \GTM(\prt\Gw)$.
\es\medskip

Since $\prt \Gw$ is $C^2$, there exists $\gd^*>0$ such that for any $\gd \in (0,\gd^*]$ and $x\in\Gw$ such that $d(x)<\gd$, there exists a unique $\gs(x)\in \prt \Gw$ such that $|x-\gs(x)|=d(x)$. We set $\gs(x)=Proj_{_{\prt \Gw}}(x)$. Furthermore, if ${\bf n}={\bf n}_{\gs(x)}$ is the normal outward unit vector to $\prt \Gw$ at $\gs(x)$, we have $x=\gs(x)-d(x){\bf n}_{\gs(x)}$. For $\gd \in (0,\gd^*]$, we set
	$$ \BA{ll} \Gw_\gd=\{x \in \Gw: d(x) \leq \gd\}, \\
						 \Gw'_\gd=\{x \in \Gw: d(x) > \gd\}, \\
						 \Gs_\gd=\prt \Gw'_\gd=\{x \in \Gw: d(x) = \gd\}, \\
						 \Gs:=\Gs_0=\prt \Gw.
	\EA $$
For any $\gd \in (0,\gd^*]$, the mapping $x\mapsto (\gd(x),\gs(x))$ defines a $C^1$ diffeomorphism from $\Gw_{\gd}$ to $(0,\gd) \ti \Gs$. Therefore we can write $x=\gs(x)-d(x){\bf n}_{\gs(x)}$ for every $x \in \Gw_\gd$. Any point $x \in \ovl \Gw_{\gd^*}$ is represented by the couple $(\gd,\gs) \in [0,\gd^*]\ti \Gs$ with formula $x=\gs-\gd \bf n_\gs$. This system of coordinates which will be made more precise in the boundary trace construction is called {\it  flow coordinates}. \medskip

\noindent{\it Proof of \rth{meas}.} 
\noindent{\it Step 1: Construction of approximate solutions.} Let $\{\gm_n\}$ be a sequence of  positive functions in $C^1(\prt \Gw)$ such that $\{\gm_n\}$ converges to $\gm$ in the weak sense of measures and $\norm{\gm_n}_{L^1(\prt \Gw)} \leq c_2\norm{\gm}_{\GTM(\prt \Gw)}$ for all $n$, where $c_2$ is a positive constant independent of $n$. We next consider the following problem
	\bel{v_n} \left\{ \BA{lll}
		- \Gd v +g(|\nabla (v + \BBP^{\Gw}[\gm_n])]) = 0 \qq &\text{ in } \Gw\\
		\phantom{- \Gd  +g(\abs{\nabla (v + \BBP^{\Gw}[\gm_n])}) }
		v = 0  &\text{ on } \prt \Gw.
	\EA \right. \ee
It is easy to see that $0$ and $-\BBP^{\Gw}[\gm_n]$ are respectively supersolution and subsolution of $(\ref{v_n})$. By \cite[Theorem 6.5]{KaKr} there exists a solution $v_n \in W^{2,p}(\Gw)$ with $1<p<\ity$ to problem $(\ref{v_n})$ satisfying $-\BBP^{\Gw}[\gm_n] \leq v_n \leq 0$. Thus the function $u_n = v_n + \BBP^{\Gw}[\gm_n]$ is a solution of
	\bel{u_n} \left\{ \BA{lll}
		- \Gd u_n + g(\abs{\nabla u_n}) = 0 \qq &\text{ in } \Gw\\
		\phantom{- \Gd  + g(\abs{\nabla u_n}) }
		u_n = \gm_n  &\text{ on } \prt \Gw.
	\EA \right. \ee
By the maximum principle, such solution is the unique solution of $(\ref{u_n})$. \medskip

\noindent{\it Step 2: We claim that $\{u_n\}$ and $\{\abs{\nabla u_n}\}$ remain uniformly bounded respectively in $ M^\frac{N}{N-1}(\Gw)$ and $M_d^\frac{N+1}{N}(\Gw)$.} 
Let $\gx$ be the solution to 
	 \bel{eta} \left\{ \BA{lll}
		- \Gd \gx = 1 \qq &\text{ in } \Gw\\
		\phantom{- \Gd}
		\gx = 0  &\text{ on } \prt \Gw,
	\EA \right. \ee
then there exists a constant $c_3>0$ such that 
	\bel{eta-est} \frac{1}{c_3}<-\frac{\prt \gx}{\prt \bf n} < c_3 \text{ and } \frac{d(x)}{c_3} \leq \gx \leq c_3d(x). \ee
By multiplying the equation in $(\ref{u_n})$ by $\gx$ and integrating on $\Gw$, we obtain
	$$ \myint{\Gw}{}u_n dx + \myint{\Gw}{}g(\abs{\nabla u_n}) \gx dx = - \myint{\prt \Gw}{}\mu_n \myfrac{\prt \gx}{\prt \bf n} dS, $$
which implies
	\bel{EM2-1} \myint{\Gw}{}u_n dx + \myint{\Gw}{}d(x) g(\abs{\nabla u_n}) dx \leq c_4\norm{\gm}_{\GTM(\prt \Gw)} \ee
where $c_4$ is a positive constant independent of $n$. By \rprop{E1} and by noticing that $u_n \leq \BBP^{\Gw}[\gm_n]$, we get 
	\bel{EM2-2} \norm{u_n}_{M^\frac{N}{N-1}(\Gw)} \leq \norm{\BBP^{\Gw}[\gm_n]}_{M^\frac{N}{N-1}(\Gw)} \leq c_1\norm{\gm_n}_{L^1(\prt \Gw)} \leq c_1c_2\norm{\gm}_{\GTM(\prt \Gw)}. \ee
Set $f_n= - g(\abs{\nabla u_n})$ then $f_n \in L_d^1(\Gw)$ and $u_n$ satisfies
	\bel{u_n*} \myint{\Gw}{}(- u_n \Gd \zeta - f_n \zeta)dx = - \myint{\prt \Gw}{}\gm_n \myfrac{\prt \zeta}{\prt \bf n} dS \ee
for any $\zeta \in X(\Gw)$.
From $(\ref{M2'})$ and \rprop{E1}, we derive that
	\bel{EM2-3} \norm{\nabla u_n}_{M_d^\frac{N+1}{N}(\Gw)} \leq c_1\left(\norm{f_n}_{L_d^1(\Gw)} + \norm{\gm_n}_{L^1(\prt \Gw)}\right), \ee
which, along with $(\ref{EM2-1})$, implies that
	\bel{EM2-4} \norm{\nabla u_n}_{M_d^\frac{N+1}{N}(\Gw)} \leq c_5\norm{\gm}_{\GTM(\prt \Gw)} \ee
where $c_5$ is a positive constant depending only on $\Gw$ and $N$.
Thus the claim follows from $(\ref{EM2-2})$ and $(\ref{EM2-4})$. \smallskip
	
\noindent{\it Step 3: Existence of a solution.} By standard results on elliptic equations and measure theory \cite[Cor. IV 27]{Br0}, the sequences $\{u_n\}$ and $\{\abs{\nabla u_n}\}$ are relatively compact in  $L^1_{loc}(\Gw)$. Therefore, there exist a subsequence, still denoted by $\{u_n\}$, and a function $u$ such that $\{u_n\}$ converges to $u$ in $L^1_{loc}(\Gw)$ and a.e. in $\Gw$.\smallskip

\noindent(i) The sequence $\{u_n\}$ converges to $u$ in $L^1(\Gw)$: let $E\subset\Gw$ be a Borel subset, then

\bel{EM2-5} \myint{E}{}u_ndx\leq |E|^{\frac{1}{N}}\norm{u_n}_{M^{\frac{N}{N-1}}(\Gw)}
\leq c_1c_2|E|^{\frac{1}{N}}\norm\gm_{\GTM(\prt\Gw)}.
\ee
The convergence of $\{u_n\}$ in $L^1(\Gw)$ follows by Vitali's theorem.\smallskip

\noindent(ii) The sequence $g(|\nabla u_n|)$ converges to $g(|\nabla u|)$ in $L_d^1(\Gw)$: consider again a Borel set $E\subset\Gw$, $\gl>0$ and write
$$\myint{E}{}d(x)g(|\nabla u_n|)dx\leq\myint{E\cap\{x:|\nabla u_n(x)|\leq\gl\}}{}d(x)g(|\nabla u_n|)dx+\myint{\{x:|\nabla u_n(x)|>\gl\}}{}d(x)g(|\nabla u_n|)dx.
$$
First
\bel{EM2-5'}\BA {l}\myint{E\cap\{x:|\nabla u_n(x)|\leq\gl\}}{}d(x)g(|\nabla u_n|)dx
\leq g(\gl) |E|_d.
\EA\ee
Then
$$\myint{E\cap\{x:|\nabla u_n(x)|>\gl\}}{}d(x)g(|\nabla u_n|)dx\leq-\myint{\gl}{\infty}g(s)d\gw_n(s)
$$
where $\gw_n(s)=|\{x\in\Gw:|\nabla u_n(x)|>s\}|_d$. Using the fact that $g'\geq 0$ combined with $(\ref{Ma'})$ and $(\ref{EM2-4})$, we get
	$$\BA{l} 		
	-\myint{\gl}{t}g(s)d\gw_n(s)=g(\gl)\gw_n(\gl)-g(t)\gw_n(t)+\myint{\gl}{t}\gw_n(s)g'(s)ds\\[4mm]\phantom{-\myint{\gl}{t}d\gw(s)}\leq 	
	g(\gl)\gw_n(\gl)-g(t)\gw_n(t)+c_6\norm{\gm}^{\frac{N+1}{N}}_{\mathfrak 	
	M(\prt \Gw)}\myint{\gl}{t}s^{-\frac{N+1}{N}}g'(s)ds\\[4mm]\phantom{-\myint{\gl}{t}d\gw(s)}
	\leq \left(\gw_n(\gl)-c_6\norm{\gm}^{\frac{N+1}{N}}_{\mathfrak M(\prt \Gw)}\gl^{-\frac{N+1}{N}}\right)g(\gl)
	-\left(\gw_n(t)-c_6\norm{\gm}^{\frac{N+1}{N}}_{\mathfrak M(\prt \Gw)}t^{-\frac{N+1}{N}}\right)g(t)\\[4mm] 
	\phantom{-----------}+c_6\frac{N+1}{N}\norm{\gm}^{\frac{N+1}{N}}_{\mathfrak M(\prt \Gw)}\myint{\gl}{t}g(s)s^{-\frac{2N+1}{N}}ds.
	\EA $$
We have already used the fact that  $\gw_n(\gl)\leq c_6\norm{\gm}^{\frac{N+1}{N}}_{\mathfrak M(\prt \Gw)}\gl^{-\frac{N+1}{N}}$, and since the condition $(\ref{M})$ holds, $\liminf_{t\to\infty}t^{-\frac{N+1}{N}}g(t)=0$. Letting $t\to\infty$ we derive
\bel{EM2-6}
\myint{E\cap\{x:|\nabla u_n(x)|>\gl\}}{}d(x)g(|\nabla u_n|)dx\leq 
c_6\frac{N+1}{N}\norm{\gm}^{\frac{N+1}{N}}_{\mathfrak M(\prt \Gw)}\myint{\gl}{\infty}g(s)s^{-\frac{2N+1}{N}}ds.
\ee
For $\ge>0$ we fix $\gl$ in order that the right-hand side of $(\ref{EM2-6})$ be smaller than $\frac{\ge}{2}$. Thus, if 
$|E|_d\leq \frac{\ge}{2g(\gl)+1}$, 
we obtain
\bel{EM2-7}
\myint{E}{}d(x)g(|\nabla u_n|)dx\leq\ge.
\ee
The convergence follows again by Vitali's theorem.
Next for any $\zeta \in X(\Gw)$, we have
	\bel{EM2-8} \myint{\Gw}{}(-u_n\Gd \zeta + g(\abs{\nabla u_n})\zeta)dx=-\myint{\prt \Gw}{}\gm_n\myfrac{\prt \zeta}{\prt \bf n}dS
	\ee
By taking into account the fact that $|\zeta|\leq cd$ in $\Gw$, we can pass to the limit in each term in $(\ref{EM2-8})$ and obtain $(\ref{M2})$; so $u$ is a solution of $(\ref{M1})$. Clearly $u \in M^\frac{N}{N-1}(\Gw)$ and  $\abs{\nabla u} \in M_d^\frac{N+1}{N}(\Gw)$ from $(\ref{M2'})$ and \rprop{E1}. \medskip
	
\noindent{\it Step 4: Existence of a maximal solution.} We first notice that any solution $u$ of $(\ref{M1})$ is smaller than $\BBP^\Gw[\gm]$. Then $u\leq \BBP^\Gw[\gm]$ in $\Gw'_\gd$ and by the maximum principle $u\leq  u_\gd$ which satisfies
\bel{EM2-7*}\left\{\BA {ll}
-\Gd  u_\gd+g(|\nabla  u_\gd|)=0\qq&\text{in }\Gw'_\gd\\[2mm]\phantom{-\Gd +g(|\nabla \bar u_\gd|)}
 u_\gd=\BBP^\Gw[\gm]\qq&\text{on }\Gs_\gd.
\EA\right.\ee
As a consequence, $0<\gd<\gd'\Longrightarrow  u_\gd\leq  u_{\gd'}$ in $\Gw'_{\gd'}$ and $ u_\gd\downarrow
\bar u_\gm$ which is not zero if $\gm$ is so, since it is bounded from below by the already constructed solution $u$. 
We extend $ u_\gd$, $|\nabla u_\gd|$ and $g(|\nabla u_\gd|)$ by zero outside $\ovl \Gw'_\gd$ and still denote them by the same expressions. Let $E\subset\Gw$ be a Borel set and put $E_\gd=E\cap \Gw'_{\gd}$ then ($\ref{EM2-5}$) becomes
\bel{EM2-5+1}\BA {l} \myint{E_\gd}{}u_\gd dx\leq |E_\gd|^{\frac{1}{N}}\norm{u_\gd}_{M^{\frac{N}{N-1}}(\Gw'_\gd)}
\leq c_1c_2|E_\gd|^{\frac{1}{N}}\norm{\BBP^\Gw[\gm]|_{_{\Gs_\gd}}}_{L^1(\Gs_\gd)}\\[4mm]\phantom{\myint{E_\gd}{}u_\gd dx}
\leq c_1c_2c_7|E|^{\frac{1}{N}}\norm{\gm}_{\GTM(\Gs)}.
\EA\ee
Set $d_\gd(x):=\dist (x,\Gw_\gd)$ ($=(d(x)-\gd)_+$ if $x\in\Gw_{\gd^*}:=\Gw\setminus\Gw'_{\gd^*}$), we have
$$\myint{E_\gd\cap\{x:|\nabla u_\gd|>\gl\}}{}d_\gd(x)g(|\nabla u_\gd|)dx\leq-\myint{\gl}{\infty}g(s)d\gw_{\gd}(s),
$$
where $\gw_\gd(s)=|\{x\in \Gw:|\nabla u_\gd(x)|>s\}|_{d_\gd}$. Since $\norm{\BBP^{\Gw}[\gm]|_{_{\Gs_\gd}}}_{L^1(\Gs_\gd)}\leq c_7\norm{\gm}_{\mathfrak M(\Gs)}$, $(\ref{EM2-5'})$ and
$(\ref{EM2-6})$ become respectively
\bel{EM2-5'+1}\BA {l}\myint{E_\gd\cap\{x:|\nabla u_\gd(x)|\leq\gl\}}{}d_\gd(x)g(|\nabla u_\gd|)dx
\leq g(\gl) |E_\gd|_{d_\gd}.
\EA\ee
and
\bel{EM2-6+1}
\myint{E_\gd\cap\{x:|\nabla u_\gd(x)|>\gl\}}{}d_\gd(x)g(|\nabla u_\gd|)dx\leq 
c_6\frac{N+1}{N}\norm{\gm}^{\frac{N+1}{N}}_{\mathfrak M}\myint{\gl}{\infty}g(s)s^{-\frac{2N+1}{N}}ds.
\ee
Combining $(\ref{EM2-5'+1})$ and $(\ref{EM2-6+1})$ and noting that $|E_\gd|_{d_\gd}\leq |E|_d$, we obtain that for any $\ge>0$ there exists $\gl>0$, independent of $\gd$ by $(\ref{EM2-5'+1})$, such that
\bel{EM2-7**}
\myint{E_\gd}{}d_\gd(x)g(|\nabla u_\gd|)dx\leq\ge
\ee
provided $|E|_d\leq \frac{\ge}{2g(\gl)+1}$. 

Finally, if $\gz\in X(\Gw)$ we denote by $\gz_\gd$ the solution of 
\bel{EM2-9}\left\{\BA {ll}
-\Gd \gz_\gd=-\Gd\gz\qq&\text{in }\Gw'_\gd\\[2mm]\phantom{-\Gd}
\gz_\gd=0\qq&\text{on }\Gs_\gd.
\EA\right.\ee
Then
	\bel{EM2-10} \myint{\Gw'_\gd}{}(- u_\gd\Gd \zeta_\gd + g(\abs{\nabla u_\gd})\zeta_\gd)dx=-\myint{\Gs_\gd}{}\myfrac{\prt \zeta_\gd}{\prt \bf n}\BBP^\Gw[\gm]dS
	\ee
Clearly $|\gz_\gd|\leq Cd_\gd$ and  $\gz_\gd\chi_{_{\Gw'_\gd}}\to\gz$ uniformly in $\Gw$ by standard elliptic estimates. Since the right-hand side of $(\ref{EM2-10})$ converges to $-\int_{\prt\Gw}{}\frac{\prt \zeta}{\prt \bf n}d\gm$, it follows by Vitali's theorem that  $\bar u_\gm$ satisfies $(\ref{M2})$.
\medskip
	
\noindent{\it Step 5: Stability.} Consider a sequence of positive bounded measures $\{\gm_n\}$ which converges weakly to
$\gm$. By estimates $(\ref{EM2-2})$ and $(\ref{EM2-4})$, $ u_{\gm_n}$ and $g(|\nabla u_{\gm_n}|)$ are relatively compact in $L_{loc}^1(\Gw)$ and respectively uniformly integrable in $L^1(\Gw)$ and $L_d^1(\Gw)$. Up to a subsequence, they converge a.e. respectively to $u$ and $g(|\nabla u|)$ for some function $u$. As in Step 3, $u$ is a solution of $(\ref{M1})$.
\qeda \medskip

 A variant of the stability statement is the following result which will be very useful in the analysis of the boundary trace. The proof is similar as Step 4 in the proof of \rth{meas}.
 \bcor{meas-intern} Let $g$ in $\CG_0$ satisfy $(\ref{M})$. Assume $\{\gd_n\}$ is a sequence decreasing to $0$ and $\{\gm_n\}$ is a sequence of positive bounded measures on $\Gs_{\gd_n}=\prt\Gw'_{\gd_n}$ which converges to $\gm$ in the weak sense of measures and let $u_{{\gm_n}}$ be solutions of $(\ref{M1})$ with boundary data $\gm_n$. Then there exists a subsequence $\{u_{{\gm_{n_k}}}\}$ of solutions of  $(\ref{M1})$ with boundary data $\gm_{n_k}$ which converges to a solution $u_\gm$ with boundary data $\gm$. 
 \es
\subsection{Boundary trace}
The construction of the boundary trace of positive solutions of $(\ref{B*})$ is a combination of tools developed in \cite{MV1}--\cite{MV3} with the help of a geometric construction from \cite{BaMa}. 

\bdef{tra} Let $\gm_\gd \in \GTM(\Gs_\gd)$ for all $\gd \in (0,\gd^*)$ and $\gm \in \GTM(\Gs)$. We say that $\gm_\gd\to \gm$ as $\gd \to 0$ in the sense of weak convergence of measures if 
\bel{E2}
\lim_{\gd\to 0}\myint{\Gs_\gd}{}\gf(\gs(x))d\gm_\gd=\myint{\Gs}{}\gf d\gm\qq\forall\gf\in C_c(\Gs).
\ee
A function $u\in C(\Gw)$ possesses a measure boundary trace $\gm\in \mathfrak M(\Gs)$ if
\bel{E3}
\lim_{\gd\to 0}\myint{\Gs_\gd}{}\gf(\gs(x))u(x)dS=\myint{\Gs}{}\gf d\gm\qq\forall\gf\in C_c(\Gs).
\ee
Similarly, if $A$ is a relatively open subset of $\Gs$, we say that $u$ possesses a trace $\gm$ on $A$ in the sense of weak convergence of measures if $\gm \in \GTM(A)$ and $(\ref{E3})$ holds for every $\gf \in C_c(A)$.
\es

We recall the following result \cite[Cor 2.3]{MV4}, adapted here to $(\ref{B*})$, 
\bprop{tra-reg} Assume $g:\BBR_+\to \BBR_+$ and let $u\in C^2(\Gw)$ be a positive solution of $(\ref{B*})$. Suppose that for some $z\in\prt\Gw$ there exists an open neighborhood $U$ such that
\bel{E4}
\myint{U\cap \Gw}{}g(|\nabla u|)d(x) dx<\infty.
\ee
Then $u\in L^1(K\cap \Gw)$ for every compact set $K\sbs U$ and there exists a positive Radon measure $\gn$ on $\Gs\cap U$ such that 
\bel{E5}
\lim_{\gd\to 0}\myint{\Gs_\gd\cap U}{}\gf(\gs(x))u(x)dS=\myint{\Gs \cap U}{}\gf d\gn\qq\forall\gf\in C_c(\Gs\cap U).
\ee
\es

\bdef{reg-sing} Let $u\in C^2(\Gw)$ be a positive solution of $(\ref{B*})$. A point $z\in\prt\Gw$ is a regular boundary point of $u$ if there exists an open neighborhood $U$ of $z$ such that $(\ref{E4})$ holds. The set of regular points is denoted by $\CR(u)$. Its complement $\CS(u)=\prt\Gw\setminus\CR(u)$ is called the singular boundary set of $u$.
\es

Clearly $\CR(u)$ is relatively open and there exists a positive Radon measure $\gm$ on $\CR(u)$ such that $u$ admits $\gm:=\gm(u)$ as a measure boundary trace on $\CR(u)$ and $\gm(u)$ is uniquely determined. The couple $(\CS(u),\gm)$ is called the {\it boundary trace} of $u$ and denoted by $tr_{\prt\Gw}(u)$. \medskip

The main question is to determine the behaviour of $u$ near $\CS(u)$. The following result is proved in \cite[Lemma 2.8]{MV4}.
\bprop{L1} Assume $g:\BBR_+\to \BBR_+$ and $u\in C^2(\Gw)$ be a positive solution of $(\ref{B*})$ with the singular boundary set $\CS(u)$. If  $z\in  
\CS(u)$ is such that there exists an open neighborhood $U'$ of $z$ such that $u\in L^1(U'\cap\Gw)$, then for every neighborhood $U$ of $z$ there holds
\bel{E6}
\lim_{\gd\to 0}\myint{\Gs_\gd\cap U}{}u(x)dS=\infty.
\ee
\es
\bcor{3/2} Let $u\in C^2(\Gw)$ is a positive solution of $(\ref{B})$ with $\frac{3}{2}<q\leq 2$. Then $(\ref{E6})$ holds for every $z\in \CS(u) $.
\es
\Proof This is a direct consequence of \rlemma{estfunct} since $\frac{q-2}{q-1}>-1$ implies $u\in L^1(\Gw)$.\qeda\medskip

We prove below that this result holds for any $1<q\leq 2$. 
\bth{gen} Assume $g:\BBR_+\to\BBR_+$ is continuous and satisfies
\bel{E6+1}
\liminf_{r\to \infty}\frac{g(r)}{r^q}>0
\ee
where $1<q\leq 2$. If $u\in C^2(\Gw)$ is a positive solution of $(\ref{B*})$, then $(\ref{E6})$ holds for every $z\in \CS(u) $.
\es
\Proof Up to rescaling we can assume that $g(r)\geq r^q-\gt$ for some $\gt\geq 0$. 
We recall some results from \cite{BF} in the form exposed in \cite[Sect 2]{BaMa}. There exist an open cover $\{\Gs_j\}_{j=1}^k$  of $\Gs$, an open set $\CD$ of $\BBR^{N-1}$ and  $C^2$ mappings $T_j$ from $\CD$ to  $\Gs_j$ with rank $N-1$ such that for each $\gs\in\Gs_j$ there exists a unique $a\in\CD$ with the property that $\gs=T_j(a)$. The couples $\{\CD,T^{-1}_j\}$ form a system of local charts of $\Gs$. If we set $\Gw_j=\{x\in\Gw_{\gd^*}:\gs(x)\in\Gs_j\}$ then for any $j=1,...,k$ the mapping 
$$\Gp_j: (\gd, a)\mapsto x=T_j(a)-\gd{\bf n}$$
where ${\bf n}$ is the outward unit normal vector to $\Gs$ at $T_j(a)=\gs(x)$ is a $C^2$ diffeomorphism from 
$(0,\gd^*)\ti \CD$ to $\Gw_j$. The Laplacian obtains the following expressions in terms of this system of flow coordinates provided the lines $\gs_i=ct$ are the vector fields of the principal curvatures $\bar\gk_i$ on $\Gs$
\bel{E7}
\Gd=\Gd_\gd+\Gd_\gs
\ee
where
\bel{E8}
\Gd_\gd=\frac{\prt^2}{\prt\gd^2}-(N-1)H\frac{\prt }{\prt\gd}
\ee
with $H=H(\gd,.)=\frac{1}{N-1}\sum_{i=1}^{N-1}\frac{\bar\gk_i}{1-\gd\bar\gk_i}$ being the mean curvature of $\Gs_\gd$ and 
\bel{E9}
\Gd_\gs=\frac{1}{\sqrt {|\Gl|}}\sum_{i=1}^{N-1}\frac{\prt}{\prt\gs_i}\left(\frac{\sqrt {|\Gl|}}{\bar \Gl_{ii}(1-\gd\bar\gk_i+\gk_{ii}\gd^2)}\frac{\prt}{\prt\gs_i}\right).
\ee
In this expression, $\bar \Gl= (\bar \Gl_{ij})$ is the metric tensor on $\Gs$ and it is diagonal by the choice of coordinates and 
$|\Gl|=\Gp_{i=1}^{N-1}\bar \Gl_{ii}(1-\gd\bar\gk_i)^2$. In particular 
\bel{E10}
|\nabla \xi|^2=\sum_{i=1}^{N-1}\frac{\xi_{\gs_i}^2}{\bar \Gl_{ii}(1-\gd\bar\gk_i+\gk_{ii}\gd^2)}+\xi_\gd^2
\ee
and
\bel{E11}
\nabla \xi.\nabla\eta=\sum_{i=1}^{N-1}\frac{\xi_{\gs_i}\eta_{\gs_i}}{\bar \Gl_{ii}(1-\gd\bar\gk_i+\gk_{ii}\gd^2)}+\xi_\gd\eta_\gd
=\nabla_\gs\xi.\nabla_\gs\eta+\xi_\gd\eta_\gd.\ee

 If $z\in \CS(u) $ we can assume that $U_\Gs:=U\cap \Gs$ is smooth and contained in a single chart $\Gs_j$. Let $\gf$ be the first eigenfunction of $\Gd_\gs$ in $W^{1,2}_0(U_\Gs)$ normalized so that $\max_{_{U_\Gs}}\gf=1$ and $\ga>1$ to be made precise later on. From  
 $-\Gd_\gd u-\Gd_\gs u+\frac{1}{2}(|\nabla u|^q-\gt) + \frac{1}{2}g(\abs{\nabla u})\leq 0 $, we obtain
by multiplying by $\gf^\ga$ and integrating over  $U_\Gs$
\bel{E12}\BA {l}\displaystyle-\frac{d^2}{d\gd^2}\myint{U_\Gs}{}u\gf^\ga dS+(N-1)\myint{U_\Gs}{}\frac{\prt u}{\prt\gd}\gf^\ga H dS + \ga\myint{U_\Gs}{}\gf^{\ga-1}\nabla_\gs u .\nabla_\gs\gf \, dS
\\[4mm]\phantom{---------}
+\myfrac{1}{2}\myint{U_\Gs}{}\gf^{\ga}(|\nabla u|^{q}-\gt)dS + \myfrac{1}{2}\myint{U_\Gs}{}\gf^{\ga}g(|\nabla u|)dS\leq 0.
\EA\ee
Provided $\ga>q'-1$ we obtain by H\"older inequality
\bel{E12+}\BA {l}\displaystyle\abs{\myint{U_\Gs}{}\gf^{\ga-1}\nabla_\gs u .\nabla_\gs\gf dS}
\leq \left(\myint{U_\Gs}{}|\nabla u|^q\gf^\ga dS\right)^{\frac{1}{q}}\left(\myint{U_\Gs}{}|\nabla_\gs \gf|^{q'}\gf^{\ga-q'} dS\right)^{\frac{1}{q'}}\\[4mm]\displaystyle\phantom{\displaystyle\abs{\myint{U_\Gs}{}\gf^{\ga-1}\nabla_\gs u .\nabla_\gs\gf dS}}
\leq \ge\myint{U_\Gs}{}|\nabla u|^q\gf^\ga dS+{\ge^{\frac{1}{1-q}}}\myint{U_\Gs}{}|\nabla_\gs \gf|^{q'}\gf^{\ga-q'} dS,
\EA\ee
and
\bel{E13}\BA {l}\displaystyle
\abs{\myint{U_\Gs}{}\frac{\prt u}{\prt\gd}\gf^\ga H dS}\leq  {\ge }\norm{H}_{L^\ity}\myint{U_\Gs}{}|\nabla u|^q\gf^\ga dS+{\ge^{\frac{1}{1-q}}}\norm H_{L^\ity}\myint{U_\Gs}{}\gf^{\ga} dS
\EA\ee
with $\ge>0$. We derive, with $\ge$ small enough,
\bel{E14}\BA {l}\displaystyle\frac{d^2}{d\gd^2}\myint{U_\Gs}{}u\gf^\ga dS
\geq \left(\myfrac{1}{2}-c_8\ge\right)\myint{U_\Gs}{}|\nabla u|^q\gf^\ga dS + \myfrac{1}{2}\myint{U_\Gs}{}\gf^{\ga}g(|\nabla u|)dS-c_8'
\EA\ee
where $c_8=c_8(q,H)$ and $c'_8=c'_8(N,q,H)$. Integrating $(\ref{E14})$ twice yields to
\bel{E15}\BA {l}\displaystyle\myint{U_\Gs}{}u(\gd,.)\gf^\ga dS
\geq \left(\myfrac{1}{2}-c_8\ge\right)\myint{\gd}{\gd^*}\myint{U_\Gs}{}|\nabla u|^q\gf^\ga dS(\gt-\gd)d\gt + \myfrac{1}{2}\myint{U_\Gs}{}\gf^{\ga}g(|\nabla u|)dS - c_8''.
\EA\ee
Since $z\in\CS(u)$, the right-hand side of $(\ref{E15})$ tends monotically to $\infty$ as $\gd\to 0$, which implies
that $(\ref{E6})$ holds.
\qeda\medskip

\noindent\Remark It is often usefull to consider the couple $(\CS(u),\gm)$ defining the boundary trace of $u$ as an outer regular Borel measure $\gn$ uniquely determined by
\bel{E16}\gn(E)=\left\{\BA {ll}\gm(E)\qq&\text{ if }E\subset\CR(u)\\
\infty\qq&\text{ if }E\cap\CS(u)\neq\emptyset
\EA\right.\ee
for all Borel set $E\subset\prt\Gw$, and we will denote $tr_{\prt\Gw}(u)=\gn(u)$.\medskip

The integral blow-up estimate $(\ref{E6})$ remains valid if $g \in \CG_0$ and the growth estimate $(\ref{E6+1})$ is replaced by $(\ref{M})$.
\bth{gen+M} Assume $g\in\CG_0$ satisfies $(\ref{M})$.
 If $u\in C^2(\Gw)$ is a positive solution of $(\ref{B*})$, then $(\ref{E6})$ holds for every $z\in \CS(u)$.\es
 \Proof By translation we assume $z=0\in\CS(u)$ and $(\ref{E6})$ does not hold. We proceed by contradiction, assuming that there exists an open neighborhood $U$ of $z$ such that
  \bel{E17}
  \liminf_{\gd\to 0}\myint{\Gs_\gd\cap U}{}udS<\infty.
\ee
 By \rprop{L1}, for any neighborhood $U'$ of $z$ there holds
   \bel{E18}
  \myint{\Gw\cap U'}{}udx=\infty,
\ee
which implies
 \bel{E19}
  \limsup_{\gd\to 0}\myint{\Gs_\gd\cap U'}{}udS=\infty.
\ee
For $n\in\BBN_*$, we take $U'=B_{\frac{1}{n}}$; there exists a sequence $\{\gd_{n,k}\}_{k\in\BBN}$ satisfying
$\lim_{k\to \infty}\gd_{n,k}=0$ such that 
 \bel{E19*}
  \lim_{k\to \infty}\myint{\Gs_{\gd_{n,k}}\cap B_{\frac{1}{n}}}{}udS=\infty.
\ee
Then, for any $\ell>0$, there exists $k_\ell:=k_{n,\ell}\in\BBN$ such that 
 \bel{E20}
 k\geq k_\ell\Longrightarrow \myint{\Gs_{\gd_{n,k}}\cap B_{\frac{1}{n}}}{}udS\geq\ell
 \ee
and $k_{n,\ell}\to\infty$ when $n\to \infty$.
 In particular there exists $m:=m(\ell,n)>0$ such that
 \bel{E21}
 \myint{\Gs_{\gd_{n,k_\ell}}\cap B_{\frac{1}{n}}}{}\inf\{u,m\}dS=\ell.
 \ee
By the maximum principle $u$ is bounded from below in $\Gw'_{\gd_{n,k_\ell}}$ by the solution $v:=v_{\gd_{n,k_\ell}}$ of
 \bel{E22}\left\{\BA {ll}
 -\Gd v+g(|\nabla v|)=0\qq&\text{in }\Gw'_{\gd_{n,k_\ell}}\\[2mm]
 \phantom{ -\Gd +g(|\nabla v|)}
v=\inf\{u,m\}\qq&\text{on }\Gs_{\gd_{n,k_\ell}}.
\EA\right. \ee
When $n\to\infty$, $\inf\{u,m(\ell,n)\}dS$ converges in the weak sense of measures to $\ell\gd_0$. By \rcor{meas-intern} there exists a solution $u_{\ell\gd_0}$ such that $v_{\gd_{n,k_\ell}}\to u_{\ell\gd_0}$ when $n\to\infty$ and consequently
$u\geq u_{\ell\gd_0}$ in $\Gw$. Even if $u_{\ell\gd_0}$ may not be unique, this implies
 \bel{E23}
\liminf_{\gd\to 0}\myint{\Gs_{\gd}}{}u\gz(x)dS\geq \lim_{\gd\to 0}\myint{\Gs_{\gd}}{}u_{\ell\gd_0}\gz(x)dS=\ell \ee
for any nonnegative $\gz\in C^{\infty}(\BBR^N)$ such that $\gz=1$ in a neighborhood of $0$. Since $\ell$ is arbitrary we obtain
 \bel{E24}
\liminf_{\gd\to 0}\myint{\Gs_{\gd}}{}u\gz(x)dS=\infty
\ee
which contradicts $(\ref{E17})$.\qeda
\mysection{Boundary singularities}


\subsection{Boundary data unbounded measures}
Since the works of Keller \cite{Ke} and Osserman \cite{Os}, universal a priori estimates became classical in the study of nonlinear elliptic equations with a superlinear absorption. Similar results holds for posiitive solutions of $(\ref{B})$ under some restrictions. We recall that for any $q>1$, any  solution $u$ of $(\ref{B})$ bounded from below satisfies  \cite[Th A1]{LaLi} the following estimate: for any $\ge>0$, there exists $C_\ge>0$ such that 
\bel{C1}
\sup_{d(x)\geq \ge}|\nabla u(x)|\leq C_\ge.
\ee
Later on Lions gave in \cite[Th IV 1]{Li1} a more precise estimate that we recall below.
\blemma{estgrad} Assume $q>1$ and  $u\in C^2(\Gw)$ is any solution of $(\ref{B})$ in $\Gw$. Then 
\bel{C2}
|\nabla u(x)|\leq C_1(N,q)(d(x))^{-\frac{1}{q-1}}\qq \forall x\in \Gw.
\ee
\es
Similarly, the following result is proved in \cite{Li1}.
\blemma{estfunct} Assume $q>1$ and  $u\in C^2(\Gw)$ is a solution of $(\ref{B})$ in $\Gw$. Then
\bel{C5}
|u(x)|\leq \frac{C_2(N,q)}{2-q}\left((d(x))^{\frac{q-2}{q-1}}-{\gd^*}^{\frac{q-2}{q-1}}\right)+\max\{|u(z)|: z\in \Gs_{\gd^*}\}\qq\forall x\in\Gw 
\ee
 if $q\neq 2$, and 
\bel{C6}
|u(x)|\leq C_3(N)\left(\ln\gd^*-\ln d(x)\right)+\max\{|u(z)|: z\in \Gs_{\gd^*}\}\qq\forall x\in\Gw 
\ee
if $q=2$, for some $C_2(N,q),C_3(N)>0$.
\es
\Proof Put $M_{\gd^*}:=\max\{|u(z)|: z\in \Gs_{\gd^*}\}$ and let $x\in \Gw_{\gd^*}$, $x=\gs(x)-d(x){\bf n}_{\gs(x)}$, and $x_0=\gs(x)-\gd^*{\bf n}_{\gs(x)}$. Then, using \rlemma{estgrad} and the fact that $\gs(x)=\gs(x_0)$, 
\bel{C7}
\BA{l}\abs{u(x)}\leq M_{\gd^*}+\myint{0}{1}\abs{\frac{d}{dt}u(tx+(1-t)x_0)}dt\\[2mm]\phantom{\abs{u(x)}}
\leq M_{\gd^*}+C_1(N,q)\myint{0}{1}(td(x)+(1-t)\gd^*)^{-\frac{1}{q-1}}
(\gd^*-d(x))dt.
\EA\ee
Thus we obtain $(\ref{C5})$ or $(\ref{C6})$ according to the value of $q$.\qeda\medskip

If $q=2$ and $u$ solves $(\ref{B})$, $v=e^u$ is harmonic and positive while if $q>2$, any solution remains bounded in $\Gw$. Although this last case is interesting in itself, we will consider only the case $1<q<2$.

\blemma{estfunct2} Assume $1<q<2$, $0\in\prt\Gw$ and  $u\in C(\overline \Gw\setminus \{0\})\cap C^2(\Gw)$ is a solution of $(\ref{B})$ in $\Gw$ which vanishes on $\prt\Gw\setminus \{0\}$. Then
\bel{C8}
u(x)\leq C_4(q)|x|^{\frac{q-2}{q-1}}\qq\forall x\in\Gw.
\ee

\es
\Proof For $\ge>0$, we set
$$P_\ge(r)=\left\{\BA{ll}
0&\text{if }r\leq\ge\\
\frac{-r^4}{2\ge^3}+\frac{3r^3}{\ge^2}-\frac{6r^2}{\ge}+5r-\frac{3\ge}{2}\q&\text{if }\ge<r<2\ge\\
r-\frac{3\ge}{2}&\text{if } r \geq 2\ge
\EA\right.
$$
and let $u_\ge$ be the extension of $P_\ge(u)$ by zero outside $\Gw$. There exists $R_0$ such that $\Gw\subset B_{R_0}$. Since $0\leq P'_\ge(r)\leq 1$ and $P_\ge$ is convex, $u_\ge \in C^2(\BBR^N)$ and it satisfies $-\Gd u_\ge+|\nabla u_\ge|^q\leq 0$. Furthermore $u_\ge$ vanishes in $ B^c_{R_0}$. For $R\geq R_0$ we set
$$U_{\ge,R}(x)=C_4(q)\left((|x|-\ge)^{\frac{q-2}{q-1}}-(R-\ge)^{\frac{q-2}{q-1}}\right)\qq\forall x\in B_R\setminus B_\ge,
$$
where $C_4(q)=(q-1)^{\frac{q-2}{q-1}}(2-q)^{-1}$, then $-\Gd U_{\ge,R}+|\nabla U_{\ge,R}|^q\geq 0$. Since $u_\ge$ vanishes on $\prt B_R$ and is finite on $\prt B_\ge$ it follows $u_\ge\leq U_{\ge,R}$ in $B_R \sms {\ovl B_\ge}$. Letting successively $\ge\to 0$ and $R\to\infty$ yields to $(\ref{C8})$.\qeda\medskip

Using regularity we can improve this estimate
\blemma{estfunct3} Under the assumptions of \rlemma{estfunct2} there holds
\bel{C9}
|\nabla u(x)|\leq C_5(q,\Gw)|x|^{-\frac{1}{q-1}}\qq\forall x\in\Gw.
\ee
and
\bel{C11}
u(x)\leq C_6(q,\Gw)d(x)|x|^{-\frac{1}{q-1}}\qq\forall x\in\Gw.
\ee
\es
\Proof For $\ell>0$, we set 
	\bel{scal} 
	T_{\ell}[u](x)=\ell^\frac{2-q}{q-1}u(\ell x) \forevery x \in \Gw^{\ell}:=\frac{1}{\ell}\Gw. 
	\ee
If $x\in\Gw$, we set $|x|=d$ and $u_d(y)=T_d[u](y)=d^{\frac{2-q}{q-1}}u(dy)$. Then $u_d$ satisfies $(\ref{B})$ in $\Gw^{d}=\frac{1}{d}\Gw$. Since $d\leq d^*:=\rm{diam}(\Gw)$, the curvature of $\prt\Gw^d$ is uniformly bounded and therefore standard a priori estimates (see e.g. \cite{GT}) imply that there exists $c$ depending on the curvature of $\Gw^d$ and $\max\{|u_d(y)|:\frac{1}{2}\leq |y|\leq \frac{3}{2} \}$ such that 
\bel{C10}
|\nabla u_d(z)|\leq c\qq\forall z\in\Gw^d,\frac{3}{4}\leq |z|\leq \frac{5}{4}.
\ee
By $(\ref{C8})$, $c$ is uniformly bounded. Therefore $|\nabla u (dz)|\leq cd^{-\frac{1}{q-1}}$ which implies $(\ref{C9})$. Finally, $(\ref{C11})$ follows from $(\ref{C8})$ and $(\ref{C9})$.
\qeda\medskip

In the next statement we obtain a local estimate of positive solutions which vanish only on a part of the boundary.

\bprop{local} Assume $1<q<2$. Then there exist $0<r^*\leq\gd^* $ and $C_7>0$ depending on $N$, $q$ and $\Gw$ such that for compact set $K\subset \prt\Gw$, $K\neq\prt\Gw$ and any positive solution  $u\in C(\ovl \Gw\setminus K)\cap C^2(\Gw)$ vanishing on $\prt\Gw\setminus K$ of $(\ref{B})$, there holds
\bel{L1}
u(x)\leq C_7d(x)(d_K(x))^{-\frac{1}{q-1}}\qq\forall x\in \Gw\txt{ s.t. } d(x)\leq r^*,
\ee
where $d_K(x)=\dist(x,K)$.
\es
\Proof The proof is based upon the construction of local barriers in spherical shells. We fix $x\in \Gw$ such that $d(x)\leq \gd^*$ and $\gs(x):={\rm Proj}_{\prt\Gw}(x)\in \prt\Gw\setminus K$. Set $r=d_K(x)$ and consider $\frac{3}{4}r<r'<\frac{7}{8}r$, $\gt\leq 2^{-1}r'$ and $\gw_x=\gs(x)+ \gt{\bf n}_x$. Since $\prt\Gw$ is $C^2$, there exists $r^*\leq \gd^*$, depending only on $\Gw$ such that $d_K(\gw_x)>\frac{7}{8}r$ provided $d(x)\leq r^*$.
For $A,B>0$ we  define the functions
$s\mapsto\tilde v(s)=A(r'-s)^{\frac{q-2}{q-1}}-B$ and $y\mapsto v(y)=\tilde v(|y-\gw_x|)$ respectively in $[0,r')$ and $B_{r'}(\gw_x)$. Then
$$ \BA {ll}-\tilde v''(s)-\myfrac{N-1}{s}\tilde v'(s)+|\tilde v'(s)|^q \\
\phantom{-------} 
=A\myfrac{2-q}{q-1}(r'-s)^{-\frac{q}{q-1}}\left(-\myfrac{1}{q-1}-\myfrac{(N-1)(r'-s)}{s} + \left(\myfrac{(2-q)A}{q-1}\right)^{q-1}\right).
\EA $$
We choose $A$ and $\gt>0$ such that
\bel{L2}
\frac{1}{q-1}-1+N+\frac{(N-1)r'}{\gt}\leq\left(\frac{(2-q)A}{q-1}\right)^{q-1}
\ee
so that inequality $-\Gd v+|\nabla v|^q\geq 0$ holds in $B_{r'}(\gw_x)\setminus B_\gt(\gw_x)$. We choose $B$ so that $v(\gs(x))=\tilde v(\gt)=0$, i.e. $B=A(r'-\gt)^{\frac{q-2}{q-1}}$. Since $\gt\leq\gd^*$, $B_\gt(\gw_x)\subset \Gw^c$ therefore $v\geq 0$ on $\prt\Gw \cap B_{r'}(\gw_x)$ and $v\geq u$ on $\Gw\cap \prt B_{r'}(\gw_x)$. By the maximum principle we obtain that $u\leq v$ in $\Gw\cap B_{r'}(\gw_x)$ and  in particular $u(x)\leq v(x)$ i.e. 
\bel{L3}
u(x)\leq A\left((r'-\gt-d(x))^{\frac{q-2}{q-1}}-(r'-\gt)^{\frac{q-2}{q-1}}\right)\leq \frac{A(2-q)}{q-1}(r'-\gt-d(x))^{-\frac{1}{q-1}}d(x).
\ee
If we take in particular $\gt=\frac{r'}{2}$ and $d(x)\leq \frac{r}{4}$, then $A=A(N,q)$ and
\bel{L4}
u(x)\leq c_9r'^{-\frac{1}{q-1}}d(x).
\ee
where $c_9=c_9(N,q)$. If we let $r'\to \frac{7}{8}r$ we derive $(\ref{L1})$. Next, if $x\in \Gw$ is such that $d(x)\leq \gd^*$ and $d(x) > \frac{1}{4}d_K(x)$, we combine $(\ref{L1})$ with Harnack inequality \cite{Tr2}, and a standard connectedness argument we obtain that $u(x)$ remains locally bounded in $\Gw$, and the bound on a compact subset $G$ of $\Gw$ depends only on $K$, $G$, $N$ and $q$. Since $d_K(x) \geq d(x) > \frac{1}{4}d_K(x)$ it follows from \rlemma{estfunct} that $(\ref{L1})$ holds. Finally $(\ref{L1})$ holds for every $x \in \Gw$ satisfying $d(x) \leq r^*$.
\qeda\medskip

As a consequence we have existence of positive solutions of $(\ref{B})$ in $\Gw$ with a locally unbounded boundary trace.

\bcor{localcor} Assume $1<q<q_c$. Then for any compact set $K\subsetneq \prt\Gw$, there exists a positive solution $u$ of $(\ref{B})$ in $\Gw$ such that $tr_{\prt\Gw}(u)=(\CS(u),\gm(u))=(K,0)$.
\es
\Proof For any $0<\ge$, we set $K_\ge=\{x\in\prt\Gw:d_K(x)< \ge\}$ and let $\psi_\ge$ be a sequence of smooth functions defined on $\prt\Gw$ such that $0\leq \psi_\ge\leq 1$, $\psi_\ge=1$ on $K_{\ge}$, $\psi_\ge=0$ on $\prt\Gw\setminus K_{2\ge}$ ($\ge<\ge_0$ so that $\prt\Gw\setminus K_{2\ge}\neq\emptyset$). Furthermore we assume that $\ge<\ge'<\ge_0$ implies $\psi_\ge\leq\psi_{\ge'}$. For $k\in\BBN^*$ let $u=u_{k,\ge}$ be the solution of 
\bel{L5}\left\{\BA {ll}
-\Gd u+|\nabla u|^q=0\qq&\text{in }\,\Gw\\[2mm]
\phantom{-\Gd +|\nabla u|^q}
u=k\psi_\ge\qq&\text{on }\,\prt\Gw.
\EA\right.\ee
By the maximum principle $(k,\ge)\mapsto u_{k,\ge}$ is increasing. Combining \rprop{local} with the same Harnack inequality argument as above we obtain that $u_{k,\ge}(x)$ remains locally bounded in $\Gw$ and satisfies $(\ref{L1})$, independently of $k$ and $\ge$. By regularity it remains locally compact in the $C^1$-topology of $\overline {\Gw}\setminus K$. If we set $u_{\infty,\ge}=\lim_{k\to\infty}u_{k,\ge}$, then it is a solution of $(\ref{B})$ in $\Gw$ which satisfies 
$$\lim_{x\to y\in K_\ge}u_{\infty,\ge}(x)=\infty\qq\forall\,y\in K_\ge,
$$
locally uniformly in $K_\ge$. Furthermore, if $y\in K_\ge$ is such that $\overline{B_{\gth}(y)}\cap\prt\Gw\subset K_\ge$ for some $\gth>0$, then for any $k$ large enough there exists $\gth_k<\gth$ such that 
$$\myint{\prt\Gw}{}\chi_{_{\overline{B_{\gth_k}(y)}\cap\prt\Gw}}dS=k^{-1}. 
$$
For any $\ell>0$, $u_{k\ell,\ge}$ is bounded from below by $u:=u_{k\ell,{B_{\gth_k}(y)}\cap\prt\Gw}$ which satisfies
\bel{L5*}\left\{\BA {ll}
-\Gd u+|\nabla u|^q=0\qq&\text{in }\,\Gw\\[2mm]
\phantom{-\Gd +|\nabla u|^q}
u=k\ell\chi_{_{\overline{B_{\gth_k}(y)}\cap\prt\Gw}}\qq&\text{on }\,\prt\Gw.
\EA\right.\ee
When $k\to\infty$, $u_{k\ell,B_{\gth_k}(y)}$ converges to $u_{\ell\gd_y}$ by \rth{meas} for the stability and \rth{uni1} for the uniqueness. It follows that $u_{\infty,\ge}\geq u_{\ell\gd_y}$. Letting $\ge\to 0$ and using the same local regularity-compactness argument we obtain that $u_K:=u_{\infty,0}=\lim_{\ge\to 0}u_{\infty,\ge}$ is a positive solution of $(\ref{B})$ in $\Gw$ which vanishes on $\prt\Gw\setminus K$ and satisfies
$$u_K\geq u_{\ell\gd_y}\Longrightarrow\lim_{\gd\to 0}\myint{\Gs_\gd\cap B_\gt(y)}{}u_K(x)dS\geq \ell,
$$
for any $\gt>0$. Since $\gt$ and $\ell$ are arbitrary, $(\ref{E6})$ holds, which implies that $y\in \CS(u_K)$. Clearly $\gm(u_K)=0$ on $\CR(u_K)=\prt\Gw\setminus \CS(u_K)$ which ends the proof.\qeda\medskip

In the supercritical case the above result cannot be always true since there exist removable boundary compact sets (see Section 4).  The following result is proved by an easy adaptation of the ideas in the proof of \rcor{localcor}.

\bcor{localcor2} Assume $q_c\leq q<2$ and let $G\subset\prt\Gw$. We assume that the boundary $\prt_{_{\prt\Gw}} G\subset\prt\Gw$ satisfies the interior boundary sphere condition relative to $\prt\Gw$ in the sense that for any $y\in \prt_{_{\prt\Gw}} G$, there exists $\ge_y>0$ and a sphere such that $B_{\ge_y}\cap\prt\Gw\subset G$ and $y\in \overline{B_{\ge_y}}$. If $\CS:=\overline G\neq \prt\Gw$ there exists a positive solution $u$ of $(\ref{B})$ with boundary trace $(\CS,0)$.
\es

\noindent\Remark It is worth noticing that the condition for the singular set to be different from all the boundary is necessary as it is shown in a recent article by Alarc\'on-Garc\'ia-Meli\'an and Quass \cite{AGMQ}. When $q_c\leq q<2$ and $\Gth\subset\prt\Gw$ it is always possible to construct a positive solution $u_\ge$ ($\ge>0$) of $(\ref{B})$ with boundary trace $({\Gth^c_\ge},0)$, where $\Gth_\ge=\{x\in\prt\Gw:d_{\Gth}(x)<\ge\}$ and the complement is relative to $\prt\Gw$. Furthermore $\ge\mapsto u_\ge$ is decreasing. If $\Gth$ has an empty interior, \rprop{local} does not apply. We conjecture that $\lim_{\ge\to 0}u_\ge$ depends on some capacity estimates on $\Gth$.\medskip

The condition that a solution vanishes outside a compact boundary set $K$ can be weakened and replaced by a local integral estimate. The next result is fundamental for existence a solution with a given general boundary trace.

\bprop{local3} Assume $1<q<2$, $U\subset\prt\Gw$ is relatively open and $\gm\in\mathfrak M(U)$ is a positive bounded Radon measure. Then for any compact set $\Gth\subset\Gw$ there exists a constant $C_8=C_8(N,q,H,\Gth,\norm\gm_{\GTM(U)})>0$ such that any positive solution $u$ of $(\ref{B})$ in $\Gw$ with boundary trace $(\CS,\gm')$ where $\CS$ is closed, $U\subset \prt\Gw\setminus\CS:=\CR$ and $\gm'$ is a positive Radon measure on 
$\CR$ such that $\gm'|_{U}=\gm$, there holds
\bel{L6}
u(x)\leq C_8\qq\forall \,x\in \Gth.
\ee
\es
\Proof We follow the notations of \rth{gen}. Since the result is local, without loss of generality we can assume that $U$ is smooth and contained in a single chart $\Gs_j$. Estimates $(\ref{E12})$-$(\ref{E15})$ are still valid under the form

\bel{L6+}\BA {l}\displaystyle\myint{U}{}u(\gd,.)\gf^\ga dS-\myint{U}{}u(\gd^*,.)\gf^\ga dS\\[4mm]
\phantom{\myint{U_\Gs}{}u(\gd,.)\gf^\ga}
\geq (1-c_{10}\ge)\myint{\gd}{\gd^*}\myint{U}{}|\nabla u|^q\gf^\ga dS(\gt-\gd)d\gt-(\gd^*-\gd)\myint{U}{}\frac{\prt u}{\prt\gd}(\gd^*,.)\gf^\ga dS-c_{10}'
\EA\ee
where $c_{10}=c_{10}(q,H)$ and $c'_{10}=c'_{10}(N,q,H)$. Since the second term in the right-hand side of $(\ref{L6+})$ is uniformly bounded by \rlemma{estgrad}, it follows that we can let 
 $\gd\to 0$ and derive, 
\bel{L7}
\myint{U}{}u(\gd^*,.)\gf^\ga dS+(1-c_{10}\ge)\myint{0}{\gd^*}\myint{U}{}|\nabla u|^q \gf^\ga \gt dSd\gt \leq \myint{U}{}\gf^\ga d\gm +c''_{10}\leq {\norm\gm}_{\GTM(U)} +c''_{10},
\ee
where $c''_{10}$ depends on the curvature $H$, $N$ and $q$. This implies that there exist some ball $B_\ga(a)$, $\ga>0$ and $a\in U$ such that $\overline{B_\ga(a)}\cap\prt\Gw\subset U$ and 
\bel{L8}
\myint{B_\ga(a)\cap \Gw}{}|\nabla u|^q d(x) dx \leq \norm\gm_{\GTM(U)} +c''_{10},
\ee
Thus, if $B_\gb(b)$ is some ball such that $\overline{B_\gb(b)}\subset B_\ga(a)\cap \Gw$, we have
\bel{L9}
\myint{B_\gb(b)}{}|\nabla u|^q dx \leq (d(b)-\gb)^{-1}\left(\norm\gm_{\GTM(U)} +c''_{10}\right).
\ee
If in  $(\ref{L6+})$  we let $\gd\to 0$ and then replace $\gd^*$ by $\gd\in (\gd_1,\gd^*]$ for $\gd_1>0$ we obtain
\bel{L11}\displaystyle\myint{U}{}\gf^\ga d\gm\geq\myint{U}{}u(\gd,.)\gf^\ga dS
 -(\gd^*-\gd)\myint{U}{}\frac{\prt u}{\prt\gd}(\gd,.)\gf^\ga dS-c'''_{10} \ee
where $c'''_{10}=c'''_{10}(N,q,H,\norm{\gm}_{\GTM(U)})$. By \rlemma{estgrad} the second term in the right-hand side remains bounded by a constant depending on $\gd_1$, $H$, $N$ and $q$. Therefore 
$\int_{U_\Gs}u(\gd,.)\gf^\ga dS$ remains bounded by a constant depending on the previous quantities and of ${\norm\gm}_{\GTM(U)}$ and consequently, assuming that $d(x)\geq \gd_1$ for all $x\in B_\gb(b)$ (i.e. $d(b)-\gb\geq \gd_1$)
\bel{L12}
u_{B_\gb(b)}:=\frac{1}{|B_\gb(b)|}\myint{B_\gb(b)}{}u dx\leq c_{11}
 \ee
where $c_{11}$ depends on $\gd_1$, $H$, $N$, $q$ and ${\norm\gm}_{\GTM(U)}$. By Poincar\'e inequality
 \bel{L13}
\left(\myint{B_\gb(b)}{}u^q dx\right)^{\frac{1}{q}}\leq c'_{11}\left[\left(\myint{B_\gb(b)}{}|\nabla u|^q dx\right)^{\frac{1}{q}}+|B_\gb(b)|^{\frac{1}{q}}u_{B_\gb(b)}\right].
 \ee
 Combining $(\ref{L9})$ and $(\ref{L12})$ we derive that $\norm u_{W^{1,q}(B_\gb(b))}$ remains bounded by a quantity depending only on $\gd_1$, $H$, $N$ and $q$ and ${\norm\gm}_{\GTM(U)}$. By the classical trace theorem in Sobolev spaces, 
 $\norm u_{L^{q}(\prt B_\gb(b))}$ remains also uniformly bounded when the above quantities are so. By the maximum principle
  \bel{L14}
u(x)\leq \BBP^{B_\gb(b)}[u|_{\prt B_\gb(b)}](x)\qq\forall\,x\in B_\gb(b),
 \ee
 where $\BBP^{B_\gb(b)}$ denotes the Poisson kernel in $B_\gb(b)$. Therefore, $u$ remains uniformly bounded in $B_{\frac{\gb}{2}}(b)$ by some constant $c''_{11}$ which also depends on ${\norm\gm}_{\GTM(U)}$, $N$, $q$, $\Gw$, $b$ and $\gb$, but not on $u$.  We end the proof by Harnack inequality and a standard connectedness argument as it has already be used in \rcor{localcor}.
 \qeda\medskip
 
 The main result of this section is the following
 \bth {borel} Assume $1<q<q_c$, $K\subsetneq \prt\Gw$ is closed and $\gm$ is a positive Radon measure on $\CR:=\prt\Gw\setminus K$. Then there exists a solution of $(\ref{B})$ such that $tr_{\prt\Gw}(u)=(K,\gm)$.
 \es
\Proof For $\ge'>\ge>0$ we set $\gn_{\ge,\ge'}=k\chi_{_{\overline K_{\ge'}}}+\chi_{_{\overline K_\ge^c}}\gm$ and denote by $u_{\ge,\ge',k,\gm}$ the maximal solution of
\bel{L15}\left\{\BA {ll}
-\Gd u+|\nabla u|^q=0\qq&\text{in }\,\Gw\\[2mm]
\phantom{-\Gd +|\nabla u|^q}
u=\gn_{\ge,\ge'}\qq&\text{on }\,\prt\Gw.
\EA\right.\ee
We recall that $K_\ge:=\{x\in\prt\Gw: d_K(x)<\ge\}$, so that $\gn_{\ge,\ge'}$ is a positive bounded Radon measure. For $0<\ge\leq \ge_0$ there exists $y\in \CR$ and $\gg>0$ such that $\overline B_\gg(y)\subset \overline K_{\ge_0}^c$. Since 
$\norm{\chi_{_{\overline K_{\ge}^c}}\gm}_{\GTM(\CR)}$ is uniformly bounded, it follows from \rprop{local3} that $u_{\ge,\ge',k,\gm}$ remains locally bounded in $\Gw$, uniformly with respect to $k$, $\ge$ and $\ge'$. Furthermore $(k,\ge,\ge')\mapsto u_{\ge,\ge',k,\gm}$ is increasing with respect to $k$. If $u_{\ge,\ge',\infty,\gm}=\lim_{k\to\infty}u_{\ge,\ge',k,\gm}$, it is a solution of $(\ref{B})$ in $\Gw$. By the same argument as the one used in the proof of \rcor{localcor}, any point $y\in K$ is such that $u_{\ge,\ge',\infty,\gm}\geq u_{\ell\gd_y}$ for any $\ell>0$. Using the maximum principle 
\bel{L16}(\ge_2\leq \ge_1,\,\ge'_1\leq \ge'_2,\,k_1\leq k_2)\Longrightarrow (u_{\ge^{}_1,\ge'_1,k^{}_1,\gm}\leq  u_{\ge^{}_2,\ge_2', k^{}_2,\gm})
\ee
Since $u_{\ge,\ge',\infty,\gm}$ remains locally bounded in $\Gw$ independently of $\ge$ and $\ge'$, we can set $u_{K,\gm}=\lim_{\ge'\to 0}\lim_{\ge\to 0}u_{\ge,\ge',\infty,\gm}$ then by the standard local regularity results $u_{K,\gm}$  is a positive solution of $(\ref{B})$ in $\Gw$. Furthermore $u_{K,\gm}>u_{\ell\gd_y}$, for any $y\in K$ and $\ell>0$; thus the set of boundary singular points of $u_{K,\gm}$ contains $K$. In order to prove that $tr_{\prt\Gw}(u_{K,\infty})=(K,\gm)$ consider a smooth relatively open set $U\subset \CR$. Using the same function $\gf^\ga$ as in \rprop{local3}, we obtain from $(\ref{L7})$ 
\bel{L17}\BA {l}\myint{U}{}u_{K,\infty}(\gd^*,.)\gf^\ga dS+(1-c_{10}\ge)\myint{0}{\gd^*}\myint{U}{}|\nabla u_{K,\infty}|^q\gf^\ga\gt dSd\gt \leq \myint{U}{}d\gm +c''_{10}.
\EA\ee
Therefore $U$ is a subset of the set of boundary regular points of $u_{K,\infty}$, which implies $tr_{\prt\Gw}(u)=(K,\gm)$ by \rprop{tra-reg}.\qeda\medskip

\noindent\Remark If $q_c\leq q<2$, it is possible to solve $(\ref{L15})$ if $\gm$ is a smooth function defined in $\CR$ and to let successively $k\to\infty$; $\ge\to 0$ and $\ge'\to 0$ using monotonicity as before. The limit function $u^*$ is a solution of $(\ref{B})$ in $\Gw$. If $tr_{\prt\Gw}(u^*)=(\CS^*,\gm^*)$, then $\CS^*\subset K$  and $\gm^*|_{\CR}=\gm$. However interior points of $K$, if any, belong to $\CS^*$ (see \rcor {localcor2}).

\subsection{Boundary Harnack inequality}
 We adapt below ideas from Bauman \cite{Ba}, Bidaut-V\'eron-Borghol-V\'eron \cite{BBV} and Trudinger \cite {Tr1}-\cite {Tr2} in order to prove a {\it boundary Harnack inequality} which is one of the main tools for analyzing the behavior of positive solutions of $(\ref{B})$ near an isolated boundary singularity. We assume that $\Gw$ is a bounded $C^2$ domain with $0\in\prt\Gw$ and $\gd^*$ has been defined for constructing the flow coordinates.

\bth{harn} Assume $0\in\prt\Gw$, $1<q<2$. Then there exist $0<r_0\leq\gd^* $ and $C_9>0$ depending on $N$, $q$ and $\Gw$ such that for any positive solution $u\in C(\Gw\cup ((\prt\Gw\sms\{0\})\cap B_{2r_0}))\cap C^2(\Gw)$ of $(\ref{B})$ vanishing on $(\prt\Gw\sms\{0\})\cap B_{2r_0}$ there holds
\bel{H1}
\frac{u(y)}{C_9d(y)}\leq \frac{u(x)}{d(x)}\leq \frac{C_9u(y)}{d(y)}
\ee
for every $x,y\in B_{\frac{2r_0}{3}}\cap \Gw$ satisfying $\frac{|y|}{2}\leq |x|\leq 2|y|$.
\es

Since $\Gw$ is a bounded $C^2$ domain, it satisfies uniform sphere condition, i.e there exists $r_0>0$ sufficiently small such that for any $x \in \prt \Gw$ the two balls $B_{r_0}(x-r_0{\bf n}_x)$ and $B_{r_0}(x+r_0{\bf n}_x)$ are subsets of $\Gw$ and $\ovl \Gw^c$ respectively. We can choose $0<r_0<\min\{\gd^*,3r^*\}$ where $r^*$ is in \rprop{local}. \medskip

We first recall the following chained property of the domain $\Gw$ \cite{Ba}.

\blemma{chain} Assume that $Q \in \prt \Gw$, $0<r<r_0$ and $h>1$ is an integer. There exists an integer $N_0$ depending only on $r_0$ such that for any points $x$ and $y$ in $\Gw \cap B_{\frac{3r}{2}}(Q)$ verifying $\min\{d(x),d(y)\} \geq r/2^h$, there exists a connected chain of balls $B_1,...,B_j$ with $j\leq N_0h$ such that 
	\bel{geo} \BA{ll} x \in B_1, y \in B_j, \q B_i\cap B_{i+1} \neq \ems  \text{ for } 1\leq i \leq j-1 \\
        \text{and } 2B_i \sbs B_{2r}(Q) \cap \Gw \text{ for } 1\leq i \leq j.
	\EA \ee
\es
The next result is an internal Harnack inequality.
\blemma{Har1} Assume $Q \in (\prt \Gw\sms\{0\}) \cap B_{\frac{2r_0}{3}}$ and $0<r\leq \abs{Q}/4$. Let $u\in C(\Gw\cup ((\prt\Gw\sms\{0\})\cap B_{2r_0}))\cap C^2(\Gw)$ be a positive solution of $(\ref{B})$ vanishing on $(\prt\Gw\sms\{0\})\cap B_{2r_0}$. Then there exists a positive constant $c_{12}>1$ depending on $N$, $q$, $\gd^*$ and $r_0$ such that 
 \bel{Har1-1} u(x) \leq c_{12}^h u(y),\ee
for every $x,y \in B_{\frac{3r}{2}}(Q)\cap \Gw$ such that $\min\{d(x),d(y)\} \geq r/2^h$ for some $h \in \BBN$.
\es
\Proof We first notice that for any $\ell>0$, $T_{\ell}[u]$ satisfies $(\ref{B})$ in $\Gw^{\ell}$ where $T_\ell$ is defined in $(\ref{scal})$. If we take in particular $\ell=|Q|$, we can assume $|Q|=1$ and the curvature of the domain $\Gw^{|Q|}$ remains bounded. By \rprop{local}  
	\bel{Har1-2} u(x) \leq C'_7\forevery x \in B_{2r}(Q)\cap \Gw \ee
where $C'_7$ depends on $N$, $q$, $\gd^*$. By \rlemma{chain} there exist an integer $N_0$ depending on $r_0$ and a connected chain of $j\leq N_0h$ balls $B_i$ with respectively radii $r_i$ and centers $x_i$, satisfying $(\ref{geo})$. Hence due to \cite[Corollary 10]{Tr1} and \cite[Theorem 1.1]{Tr2} there exists a positive constant $c'_{12}$ depending on $N$, $q$, $\gd^*$ and $r_0$ such that for every $1 \leq i \leq j$,
	\bel{Har1-3} \sup_{B_i}u \leq c'_{12} \inf_{B_i}u, \ee
which yields to $(\ref{Har1-1})$ with $c_{12}=c_{12}'^{N_0}$. \qeda

By proceeding as in \cite{Ba} and \cite{BBV}, we obtain the following results.
\blemma{Har2} Assume the assumptions on $Q$ and $u$ of \rlemma{Har1} are fulfilled. If $P \in \prt \Gw \cap B_r(Q)$ and $0<s<r$, there exist two positive constants $\gd$ and $c_{13}$ depending on $N$, $q$ and $\Gw$ such that
	\bel{Har2-1} u(x) \leq c_{13} \myfrac{\abs{x-P}^\gd}{s^\gd}M_{s,P}(u) \ee
for every $x \in B_s(P) \cap \Gw$, where $M_{s,P}(u)=\max\{u(z): z\in B_s(P)\cap \Gw\}$.
\es
\bcor{Har3} Assume $Q \in (\prt \Gw\sms\{0\}) \cap B_{\frac{2r_0}{3}}$ and $0<r\leq \abs{Q}/8$. Let $u\in C(\Gw\cup ((\prt\Gw\sms\{0\})\cap B_{2r_0}))\cap C^2(\Gw)$ positive solution of $(\ref{B})$ vanishing on $(\prt\Gw\sms\{0\})\cap B_{2r_0}$. Then there exists a constant $c_{14}$ depending only on $N$, $q$, $\gd^*$ and $r_0$ such that
	\bel{Har3-1} u(x) \leq c_{14} u(Q-\frac{r}{2}{\bf n}_{_Q}) \forevery x \in B_r(Q) \cap \Gw. \ee
\es

\blemma{LUE1} Assume $Q \in (\prt \Gw\sms\{0\}) \cap B_{\frac{2r_0}{3}}$ and $0<r\leq \abs{Q}/8$. Let $u\in C(\Gw\cup ((\prt\Gw\sms\{0\})\cap B_{2r_0}))\cap C^2(\Gw)$ positive solution of $(\ref{B})$ vanishing on $(\prt\Gw\sms\{0\})\cap B_{2r_0}$. Then there exist $a \in (0,1/2)$ and $c_{15}>0$ depending on $N$, $q$, $\gd^*$ and $r_0$ such that
	\bel{LUE1} \myfrac{1}{c_{15}}\myfrac{t}{r} \leq \myfrac{u(P-t{\bf n}_{_P})}{u(Q-\frac{r}{2}{\bf n}_{_Q})} \leq c_{15}\myfrac{t}{r}\ee
for any $P \in B_r(Q) \cap \prt \Gw$ and $0 \leq t < \frac{a}{2}r$.
\es
\medskip

\noindent{\it Proof of \rth{harn}.} Assume $x \in B_{\frac{2r_0}{3}} \cap \Gw$ and set $r=\frac{\abs{x}}{8}$. \medskip

\noindent{\it Step 1: Tangential estimate: we  suppose  $d(x) < \frac{a}{2}r$}. Let $Q \in \prt \Gw \sms \{0\}$ such that $\abs{Q}=\abs{x}$ and $x \in B_r(Q)$. By \rlemma{LUE1},
	\bel{Hark1-2} \myfrac{8}{c_{15}}\myfrac{u(Q-\frac{r}{2}{\bf n}_{_Q})}{\abs{x}} \leq \myfrac{u(x)}{d(x)} \leq 8c_{15}\myfrac{u(Q-\frac{r}{2}{\bf n}_{_Q})}{\abs{x}}. \ee
We can connect $Q-\frac{r}{2}{\bf n}_{_Q}$ with $-2r{\bf n}_{_0}$ by $m_1$ (depending only on $N$) connected balls $B_i=B(x_i,\frac{r}{4})$ with $x_i \in \Gw$ and $d(x_i) \geq \frac{r}{2}$ for every $1 \leq i \leq m_1$. It follows from $(\ref{Har1-3})$ that
	$$ c_{12}'^{-m_1}u(-2r{\bf n}_{_0}) \leq u(Q-\frac{r}{2}{\bf n}_{_Q}) \leq c_{12}'^{m_1}u(-2r{\bf n}_{_0}), $$
which, together with $(\ref{Hark1-2})$ leads to
	\bel{Hark1-3} \myfrac{8}{c_{12}'^{m_1}c_{15}}\myfrac{u(-2r{\bf n}_{_0})}{\abs{x}} \leq \myfrac{u(x)}{d(x)} \leq 8c_{12}'^{m_1}c_{15}\myfrac{u(-2r{\bf n}_{_0})}{\abs{x}}. \ee
\noindent{\it Step 2: Internal estimate: $d(x) \geq \frac{a}{2}r$}. We can connect $-2r{\bf n}_{_0}$ with $x$ by $m_2$ (depending only on $N$) connected balls $B'_i=B(x'_i,\frac{a}{4}r)$ with $x'_i \in \Gw$ and $d(x'_i) \geq \frac{a}{2}r$ for every $1 \leq i \leq m_2$. By applying again $(\ref{Har1-3})$ and keeping in mind the estimate $\frac{a}{4}\abs{x}<d(x)\leq \abs{x}$, we get 
	\bel{Hark1-4}  \myfrac{a}{4c_{12}'^{m_2}}\myfrac{u(-2r{\bf n}_{_0})}{\abs{x}} \leq \myfrac{u(x)}{d(x)} \leq \myfrac{4c_{12}'^{m_2}}{a}\myfrac{u(-2r{\bf n}_{_0})}{\abs{x}}. \ee
\noindent{\it Step 3: End of proof.} Take $\frac{\abs{x}}{2}\leq s \leq 2\abs{x}$, we can connect $-2r{\bf n}_{_Q}$ with $-s{\bf n}_{_Q}$ by $m_3$ (depending only on $N$) connected balls $B''_i=B(x''_i,\frac{r}{2})$ with $x''_i \in \Gw$ and $d(x''_i) \geq r$ for every $1 \leq i \leq m_3$. This fact, joint with $(\ref{Hark1-3})$ and $(\ref{Hark1-4})$, yields
	\bel{Hark1-5}  \myfrac{1}{C_{9}'}\myfrac{u(-s{\bf n}_{_0})}{\abs{x}} \leq \myfrac{u(x)}{d(x)} \leq C_{9}'\myfrac{u(-s{\bf n}_{_0})}{\abs{x}} \ee
where $C_{9}'=C_{9}'(N,q,\Gw)$. Finally let $y \in B_{\frac{2r_0}{3}} \cap \Gw$ satisfy $\frac{\abs{x}}{2} \leq \abs{y} \leq 2\abs{x}$. By applying twice $(\ref{Hark1-5})$ we get $(\ref{H1})$ with $C_{9}=C_{9}'^2$. \qeda\medskip

A direct consequence of \rth{harn} is the following useful form of boundary Harnack inequality.
\bcor{harn2} Let $u_i \in C(\Gw\cup ((\prt\Gw\sms\{0\})\cap B_{2r_0}))\cap C^2(\Gw)$  ($i=1,2$) be two nonnegative solutions of $(\ref{B})$ vanishing on $(\prt\Gw\sms\{0\})\cap B_{2r_0}$. Then there exists a constant $C_{10}$ depending on $N$, $q$ and $\Gw$ such that for any $r\leq \frac{2r_0}{3}$
	\bel{Hark2-1} \BA{ll}
	\sup\left(\myfrac{u_1(x)}{u_2(x)}:x \in \Gw \cap (B_r\sms B_{\frac{r}{2}})\right) \\
	\phantom{------}
	\leq C_{10}\inf \left(\myfrac{u_1(x)}{u_2(x)}:x \in \Gw \cap (B_r\sms B_{\frac{r}{2}})\right).
	\EA \ee
\es


\subsection{Isolated singularities}
\rth{meas} assert the existence of a solution to $(\ref{M1})$ for any positive Radon measure $\gm$ if $g\in\CG_0$ satisfies $(\ref{M})$, and the question of uniqueness of this problem is still an open question, nevertheless when $\gm=\gd_z$ with $z \in \prt \Gw$, we have the following result
\bth{uni1} Assume $1< q < q_c$, $z\in\prt\Gw$ and $c > 0$.  Then there exists a unique solution $u:=u_{c\gd_z}$ to
 	\bel{Pc} \left\{ \BA{lll} -\Gd u + \abs{\nabla u}^q  = 0 \qq &\text{in } \Gw \\
 	\phantom{ -\Gd  + \abs{\nabla u}^q}
 	u = c\gd_z &\text{on } \prt \Gw
 	\EA \right. \ee
	Furthermore the mapping $c\mapsto u_{c\gd_z}$ is increasing.
\es 	
\blemma{APE-grad} Under the assumption of \rth{uni1}, there holds
	\bel{APE-grad1} \abs{\nabla u(x)} \leq  C_{11}c\abs{x-z}^{-N} \forevery x \in \Gw \ee
with $C_{11}=C_{11}(N,q,\gk)>0$ where $\gk$ is the supremum of the curvature of $\prt\Gw$.
\es
\Proof Up to a translation we may assume $z=0$. By the maximum principle $0<u(x)\leq cP^\Gw(x,0)$ in $\Gw$. For $0<\ell \leq 1$, set $v_\ell=T_\ell[u]$ where $T_\ell$ is the scaling defined in $(\ref{scal})$, then $v_\ell$ satisfies
	\bel{APE-grad2} \left\{ \BA{lll}
	- \Gd v_\ell + \abs{\nabla v_\ell}^q &= 0 \qq &\text{in } \Gw^{\ell} \\
	\phantom{------}
	v_\ell &= \ell^{\frac{2-q}{q-1}+1-N}c\gd_0 &\text{on } \prt \Gw^{\ell}
	\EA \right. \ee
where $\Gw^{\ell}=\frac{1}{\ell}\Gw$ and by the maximum principle
	$$ 0 < v_\ell(x) \leq \ell^{\frac{2-q}{q-1}+1-N}cP^{\Gw^{\ell}}(x,0) \forevery x \in \Gw^{\ell}. $$
Since the curvature of $\prt \Gw^{\ell}$ remains bounded when $0<\ell\leq 1$, there holds (see \cite{Lib})
	\bel{APE-grad3} \BA{lll} \sup\{\abs{\nabla v_\ell(x)}:x\in \Gw^{\ell} \cap (B_2\sms B_{\frac{1}{2}})\} \\
	\phantom{-------}
	\leq C_{11}'\sup\{v_\ell(x): x\in \Gw^{\ell} \cap (B_3\sms B_{\frac{1}{3}})\} \\
	\phantom{-------}
	\leq C_{11}'\ell^{\frac{2-q}{q-1}}\sup\{u(\ell x): x\in \Gw^{\ell} \cap (B_3\sms B_{\frac{1}{3}})\} \\
	\phantom{-------}
	\leq C_{11}c\ell^{\frac{2-q}{q-1}+1-N}
	\EA \ee
where $C_{11}$ and $C_{11}'$ depend on $N$, $q$ and $\gk$. Consequently
	$$ \ell^{\frac{2-q}{q-1}+1}\abs{\nabla u}(\ell x) \leq C_{11}(N,q,\gk)c\ell^{\frac{2-q}{q-1}+1-N} \forevery x \in \Gw^{\ell} \cap (B_2\sms B_{\frac{1}{2}}), \q\forall \ell>0 $$
Set $\ell x = y$ and $\abs{x}=1$, then 
	$$ \abs{\nabla u(y)} \leq C_{11}\abs{y}^{-N} \forevery y \in \Gw. $$ \qeda
\blemma{lim}
	\bel{AB4} \lim_{\abs{x} \to 0}\myfrac{\BBG^{\Gw}[\abs{x}^{-Nq}]}{P(x,0)}=0. \ee
\es
We recall the following estimates for the Green fuction (\cite{BBV}, \cite{GV}, \cite{V1} and \cite{V2})
	$$ G^\Gw(x,y) \leq c_{16}d(x)\abs{x-y}^{1-N} \forevery x,y\in \Gw, x \neq y $$
and
	$$ G^\Gw(x,y) \leq c_{16}d(x)d(y)\abs{x-y}^{-N} \forevery x,y\in \Gw, x \neq y. $$
where $c_{16}=c_{16}(N,\Gw)$.
Hence, for $\ga \in (0,N+1-Nq)$, we obtain
	\bel{AB5} \BA{lll} G^\Gw(x,y) \leq \left(c_{16}d(x)\abs{x-y}^{1-N}\right)^\ga\left(c_{16}d(x)d(y)\abs{x-y}^{-N}\right)^{1-\ga} \\ 	
	\phantom{G^\Gw(x,y)}
	=c_{16}d(x)d(y)^{1-\ga}\abs{x-y}^{\ga-N} \forevery x,y\in \Gw, x \neq y,
	\EA \ee  
which follows that
	\bel{AB6} \myfrac{\BBG^{\Gw}[\abs{x}^{-Nq}]}{P^\Gw(x,0)}\leq c_{16}\abs{x}^N\myint{\BBR^N}{}\abs{x-y}^{\ga-N}\abs{y}^{1-Nq-\ga}dy \ee
By the following identity (see \cite[p. 124]{LiLo}), 
	\bel{AB7} \myint{\BBR^N}{}\abs{x-y}^{\ga-N}\abs{y}^{1-Nq-\ga}dy = c_{16}'\abs{x}^{1-Nq} \ee
where $c_{16}'=c_{16}'(N,\ga)$, we obtain
	\bel{AB8} \myfrac{\BBG^{\Gw}[\abs{x}^{-Nq}]}{P^\Gw(x,0)}\leq c_{16} c_{16}'\abs{x}^{N+1-Nq}. \ee
Since $N+1-Nq>0$, $(\ref{AB4})$ follows. \qeda \medskip

\noindent{\bf Proof of \rth{uni1}.} Since $ u=c\,\BBP^\Gw[\gd_0]-\BBG^{\Gw}[\abs{\nabla u}^q]$, 
	\bel{AB9} \lim_{\abs{x} \to 0}\myfrac{u(x)}{P^\Gw(x,0)}=c. \ee
Let $u$ and $\tl u$ be two solutions to $(\ref{Pc})$. For any $\vge>0$, set $u_\vge=(1+\vge)u$ then $u_\vge$ is a supersolution. By step 3,
	$$\lim_{x \to 0}\myfrac{u_\vge(x)}{P^\Gw(x,0)}=(1+\vge)c. $$
Therefore there exists $\gd=\gd(\ge)$ such that $u_\ge \geq \tl u$ on $\Gw \cap \prt B_\gd$. By the maximum principle, $u_\vge \geq \tl u$ in $\Gw\sms B_\gd$. Letting $\vge \to 0$ yields to $u \geq \tl u$ in $\Gw$ and the uniqueness follows. The monotonicity of $c\mapsto u_{c\gd_0}$ comes from $(\ref{AB9})$. \qeda \medskip

As a variant of the previous result we have its extension in some unbounded domains.

\bth{unbd} Assume $1<q<q_c$, and  either $\Gw=\BBR_+^N:=\{x=(x',x_N):x_N>0\}$ or $\prt\Gw$ is compact with $0\in\prt\Gw$. Then there exists one and only one solution to problem $(\ref{Pc})$.
\es
\Proof The proof needs only minor modifications in order to take into account the decay of the solutions at $\infty$. For $R>0$ we set $\Gw_R=\Gw\cap B_R$ and denote by $u:=u^R_{c\gd_0}$ the unique solution of 
 	\bel {Pc+1}\left\{ \BA{lll} -\Gd u + \abs{\nabla u}^q  = 0 \qq &\text{in } \Gw_R \\
 	\phantom{ -\Gd  + \abs{\nabla u}^q}
 	u = c\gd_0 &\text{on } \prt \Gw_R.
 	\EA \right. \ee
	Then
	 \bel {Pc+2}u^R_{c\gd_0}(x)\leq cP^{\Gw_R}(x,0)\qq\forall x\in \Gw_R. 
	 \ee	 
	Since $R\mapsto P^{\Gw_R}(.,0)$ is increasing, it follows from $(\ref{AB9})$ that  $R\mapsto u^R_{c\gd_0}$ is increasing too with limit $u^*$ and there holds 
	 	\bel {Pc+3}
	u^*(x)\leq cP^{\Gw}(x,0)\qq\forall x\in \Gw.
	\ee
	Estimate $(\ref{APE-grad1})$ is valid independently of $R$ since the curvature of $\prt\Gw_R$ is bounded (or zero if $\Gw=\BBR^N_+$). By standard local regularity theory, $\nabla u^R_{c\gd_0}$ converges locally uniformly in $\overline \Gw\setminus B_\ge$ for any $\ge>0$ when $R\to\infty$, and  thus $u^*\in C(\overline\Gw\setminus\{0\})$ is a positive solution of  $(\ref{B})$ in $\Gw$ which vanishes on $\prt\Gw\setminus\{0\}$. It admits therefore a boundary trace $tr_{\prt\Gw}(u^*)$. Estimate $(\ref{Pc+3})$ implies that $\CS(u^*)=\emptyset$ and $\gm(u^*)$ is a Dirac measure at $0$, 
which is in fact $c\gd_0$ by combining estimates $(\ref{AB9})$ for $\Gw_R$, $(\ref{Pc+2})$ and $(\ref{Pc+3})$. Uniqueness follows from the same estimate.\qeda

\medskip

We next consider the equation $(\ref{B})$ in $\BBR_+^N$. We denote by $(r,\gs) \in \BBR_+ \ti S^{N-1}$ are the spherical coordinates in $\BBR^N $ and we recall the following representation
$$S_+^{N-1}=\left\{(\sin\gf\gs',\cos\gf):\gs'\in S^{N-2},\gf\in [0,\frac{\gp}{2})\right\},
$$
$$\Gd v=v_{rr}+\frac{N-1}{r}v_r+\frac{1}{r^2}\Gd' v
$$
where $\Gd'$ is the Laplace-Beltrami operator on $S^{N-1}$,
$$\nabla v=v_r{\bf e}+\frac{1}{r}\nabla' v
$$
where $\nabla' $ denotes the covariant derivative on $S^{N-1}$ identified with the tangential derivative,
$$\Gd' v=\frac{1}{(\sin\gf)^{N-2}}\left((\sin\gf)^{N-2}v_\gf\right)_\gf+\frac{1}{(\sin\gf)^{2}}\Gd''v
$$
where $\Gd''$ is the Laplace-Beltrami operator on $S^{N-2}$. Notice that the function $\vgf_1(\gs)=\cos\gf$ is the first eigenfunction of $-\Gd'$ in $W^{1,2}_0(S^{N-1}_+)$, with corresponding eigenvalue $\gl_1=N-1$ and we choose $\gth>0$ such that $\tl \vgf_1(\gs):=\gth \cos\gf$ has mass $1$ on $S_+^{N-1}$.

We look for  a particular solution of 
	\bel{P3} \left\{ \BA{lll}
		- \Gd u + \abs{\nabla u}^q = 0 \qq &\text{ in } \BBR^N_+ \\
		\phantom{- \Gd  + \abs{\nabla u}^q}
		u = 0 &\text{ on } \prt \BBR^N_+\sms\{0\}=\BBR^{N-1}\sms \{0\}
	\EA \right. \ee
under the separable form 
	\bel{Sep-sol} u(r,\gs)=r^{-\gb}\gw(\gs) \qq (r,\gs) \in (0,\ity)\ti S_+^{N-1}.  \ee
It follows from a straightforward computation that $\gb=\frac{2-q}{q-1}$ and $\gw$ satisfies
\bel{P4} \left\{ \BA{lll}
	 \CL\gw:=-\Gd'\gw + \left((\frac{2-q}{q-1})^2 \gw^2+\abs{\nabla' \gw}^2\right)^{\frac{q}{2}}-\frac{2-q}{q-1}(\frac{q}{q-1}-N)\gw &= 0 \q &\text{in } S_+^{N-1} \\
	\phantom{------------------------,}
	\gw &= 0 &\text{on } \prt S_+^{N-1}
\EA \right. \ee
Multiplying $(\ref{P4})$ by $\vgf_1$ and integrating over $S_+^{N-1}$, we get
	$$ \BA{lll} 
	\left[N-1-\myfrac{2-q}{q-1}\bigg(\myfrac{q}{q-1}-N\bigg)\right]\myint{S_+^{N-1}}{}\gw \vgf_1 dx \\
	\phantom{---------}
	+ \myint{S_+^{N-1}}{}\left(\bigg(\myfrac{2-q}{q-1}\bigg)^2 \gw^2+\abs{\nabla' \gw}^2\right)^{\frac{q}{2}}\vgf_1 dx = 0. \EA $$
{\it Therefore  if $N-1 \geq \frac{2-q}{q-1}\big(\frac{q}{q-1}-N\big)$ and in particular if $q\geq q_c$, there exists no nontrivial solution of $(\ref{P4})$}. \medskip

In the next theorem we prove that  if $N-1 < \frac{2-q}{q-1}\big(\frac{q}{q-1}-N\big)$, or equivalently $q<\frac{N+1}{N}$, there exists a unique positive solution of $(\ref{P4})$.

\bth{exis-uniq-SS} Assume $1<q<q_c$. There exists a unique positive solution $\gw_s:=\gw \in W^{2,p}(S^{N-1}_+)$ to $(\ref{P4})$ for all $p>1$. Furthermore $\gw_s\in C^{\infty}(\overline{S^{N-1}_+})$.\es
\Proof {\it Step 1: Existence.} We first claim that $\unl \gw:=\gg_1 \vgf_1^{\gg_2}$ is a positive sub-solution of $(\ref{P4})$ where $\gg_i$ ($i=1,2$) will be determined later on. Indeed, we have
	$$ \BA{lll} 
	\CL(\unl \gw) \leq \gg_1\vgf_1^{\gg_2}\left[(N-1){\gg_2}-\myfrac{2-q}{q-1}\bigg(\myfrac{q}{q-1}-N\bigg)+2\bigg(\myfrac{2-q}{q-1}\bigg)^q\gg_1^{q-1}\vgf_1^{(q-1){\gg_2}}\right] \\[3mm]
	\phantom{\CL(\unl \gw)}  -\gg_1\vgf_1^{{\gg_2}-2}\left[\bigg(\myfrac{2-q}{q-1}\bigg)^q\gg_1^{q-1}\vgf_1^{(q-1)\gg_2+2}+\gg_2(\gg_2-1)\abs{\nabla ' \vgf_1}^2\right] 
	+\gg_1^q\gg_2^q\vgf_1^{q(\gg_2-1)}\abs{\nabla ' \vgf_1}^q \\[3mm]
	\phantom{\CL(\unl \gw)} =:\gg_1\vgf_1^{\gg_2}L_1 -  \gg_1\vgf_1^{{\gg_2}-2}L_2 + L_3.
	\EA $$
Since $q<q_c$, we can choose 
	 $$ 1< \gg_2< \myfrac{(N+q-Nq)(2-q)}{(N-1)(q-1)^2}. $$ 
Since $\vgf_1 \leq 1$, we can choose $\gg_1>0$ small enough in order that $L_1<0$ and $-\gg_1\vgf_1^{{\gg_2}-2}L_2 + L_3<0$. Thus the claim follows. 

Next, it is easy to see that $\ovl \gw=\gg_4$, with $\gg_4>0$ large enough, is a supersolution of $(\ref{P4})$ and $\ovl \gw>\unl \gw$ in $\ovl S^{N-1}_+$. Therefore there exists a solution $\gw \in W^{2,p}(S^{N-1}_+)$ to $(\ref{P4})$ such that $0< \unl \gw \leq \gw \leq \ovl \gw$ in $S_+^{N-1}$. \medskip

\noindent{\it Step 2: Uniqueness.} Suppose that $\gw_1$ and $\gw_2$ are two positive different solutions of $(\ref{P4})$ and by Hopf lemma $\nabla'\gw_i$ ($i=1,2$) does not vanish on $S^{N-1}_+$.  Up to exchanging the role of $\gw_1$ and $\gw_2$, we may assume $ \max_{S^{N-1}_+}\gw_2 \geq \max_{S^{N-1}_+}\gw_1$ and
	$$ \gl:=\inf\{c>1:c\gw_1 > \gw_2 \text{ in } S^{N-1}_+\} >  1. $$
Set $\gw_{1,\gl}:=\gl \gw_1$, then $\gw_{1,\gl}$ is a positive supersolution to problem $(\ref{P4})$. Owing to the definition of $\gw_{1,\gl}$, one of two following cases must occur. \medskip

\noindent{\it Case 1: Either $\exists \gs_0 \in S^{N-1}_+$ such that $\gw_{1,\gl}(\gs_0)=\gw_2(\gs_0)>0$ and $\nabla ' \gw_{1,\gl}(\gs_0)=\nabla ' \gw_2(\gs_0)$}. Set $\gw_{\gl}:=\gw_{1,\gl}-\gw_2$ then $\gw_{\gl}\geq 0$ in $\ovl S^{N-1}_+$, $\gw(\gs_0)=0$, $\nabla '\gw_\gl(\gs_0)=0$. Morevover,
	\bel{uniq1} \BA{ll} 
	-\Gd' \gw_\gl + (H(\gw_{1,\gl},\nabla' \gw_{1,\gl})-H(\gw_{2},\nabla' \gw_{2}))  
	 - \myfrac{2-q}{q-1}\bigg(\myfrac{q}{q-1}-N\bigg) \gw_\gl \geq 0.
	\EA \ee
where $H(s,\gx)=((\frac{2-q}{q-1})^2 s^2 + \abs{\gx}^2)^\frac{q}{2}$, $(s,\gx) \in \BBR \ti \BBR^N$. By the Mean Value theorem and $(\ref{uniq1})$, we may choose $\gg_5>0$ large enough such that
	$$ -\Gd' \gw_\gl + \myfrac{\prt H}{\prt \gx}(\ovl s, \ovl \gx)\nabla' \gw_\gl + \left[\gg_5-\myfrac{2-q}{q-1}\left(\frac{q}{q-1}-N\right)\right]\gw_\gl \geq 0 $$
where $\ovl s$ and $\ovl \gx_i$ are the functions with respect to $\gs \in S^{N-1}_+$. By the maximum principle, $\gw_\gl$ cannot achieve a non-positive minimum in $S^{N-1}_+$, which is a contradiction. \medskip

\noindent{\it Case 2: or $\gw_{1,\gl}>\gw_2$ in $S^{N-1}_+$ and $\exists \gs_0 \in \prt S^{N-1}_+$ such that}
	\bel{uniq3} \myfrac{\prt \gw_{1,\gl}}{\prt \bf n}(\gs_0)=\myfrac{\prt \gw_2}{\prt \bf n}(\gs_0). \ee
Since $\gw_{1,\gl}(\gs_0)=0$ and $\gw_{1,\gl} \in C^1(\ovl {S^{N-1}_+})$, there exists a relatively open subset $U \sbs S^{N-1}_+$ such that  $\gs_0 \in \prt U$ and 
	\bel{uniq4} \max_{\ovl U}w_{1,\gl}< q^{-\frac{1}{q-1}}\myfrac{q-1}{2-q}\left(\myfrac{q}{q-1}-N\right)^\frac{1}{q-1}.\ee
We  set $\gw_{\gl}:=\gw_{1,\gl}-\gw_2$ as in case 1. It follows that
	\bel{uniq5} \BA{ll}
	-\Gd' \gw_\gl + \myfrac{\prt H}{\prt \gx}(\ovl s, \ovl \gx)\prt_{\gs_i}\gw_\gl 
	> \myfrac{2-q}{q-1}\left[\myfrac{q}{q-1}-N - q\left(\myfrac{2-q}{q-1}\right)^{q-1}\gw_{1,\gl}^{q-1}\right]\gw_{\gl} > 0 
	\EA \ee
in $U$ owing to $(\ref{uniq4})$. By Hopf lemma $\frac{\prt \gw_\gl}{\prt \bf n}(\gs_0)<0$, which  contradicts $(\ref{uniq3})$. The regularity comes from the fact that $\gw^2+|\nabla\gw|^2>0$ in $\overline{S^{N-1}_+}$. \qeda \medskip

When $\BBR^N_+$ is replaced by a general $C^2$ bounded domain $\Gw$, the role of $\gw_s$ is crucial for describing the boundary isolated singularities. In that case we assume that $0\in\prt\Gw$ and the tangent plane to $\prt\Gw$ at $0$ is $\prt\BBR^{N-1}_+:=\{(x',0):x'\in \BBR^{N-1}\}$, with normal inward unit vector ${\bf e}_N$. If $u \in C(\ovl {\BBR^N_+}\sms\{0\})$ is a solution of $(\ref{P3})$ then so is $T_\ell[u]$ for any $\ell>0$. We say that $u$ is {\em self-similar} if $T_\ell[u] = u$ for every $\ell>0$. 

\bprop{sing} Assume $1<q<q_c$ and $0\in\prt\Gw$. Then 
\bel{uniq6} \lim_{c\to\infty}u_{c\gd_0}=u_{\infty,0}
 \ee
 where $u_{\infty,0}$ is a positive solution of $(\ref{B})$ in $\Gw$, continuous in $\ovl \Gw \sms \{0\}$ and vanishing on $\prt \Gw \sms \{0\}$. Furthermore there holds
 \bel{uniq7} 
\lim_{\tiny\BA{c}\Gw \ni x\to 0\\
\frac{x}{|x|}=\gs\in S^{N-1}_+
\EA}|x|^{\frac{2-q}{q-1}}u_{\infty,0}(x)=\gw_s(\gs),
 \ee
 locally uniformly on $S^{N-1}_+$.
\es
\Proof If $u$ is the solution of a problem $(\ref{Pc})$ in a domain $\Gth$ with boundary data  $c\gd_z$, we denote it by $u_{c\gd_z}^{\Gth}$. Let $B$ and $B'$ be two open balls tangent to $\prt\Gw$ at $0$ and such that $B\subset \Gw\subset B'^c$. Since $P^B(x,0)\leq P^\Gw(x,0)\leq P^{B'^c}(x,0)$ it follows from \rth{unbd} and $(\ref{AB9})$ that
 \bel{uniq8} 
u_{c\gd_0}^{B}\leq u_{c\gd_0}^{\Gw}\leq u_{c\gd_0}^{B'^c}.
 \ee
Because of uniqueness and whether $\Gth$ is $B$, $\Gw$ or $B'^c$, we have
 \bel{uniq9} 
T_\ell[u^\Gth_{c\gd_0}]=u^{\Gth^{\ell}}_{c\ell^{\gth}\gd_0}\qq\forall\ell>0,
 \ee
 with $\Gth^{\ell}=\frac{1}{\ell}\Gth$ and $\gth:=\frac{2-q}{q-1}+1-N$. Notice also that $c\mapsto u^\Gth_{c\gd_0}$ is increasing. Since $u^\Gth_{c\gd_0}(x) \leq C_4(q)|x|^{\frac{q-2}{q-1}}$ by $(\ref{C8})$, it follows that $u^\Gth_{c\gd_0}\uparrow u^\Gth_{\infty,0}$. As in the previous constructions, $u^\Gth_{\infty,0}$ is a positive solution of $(\ref{B})$ in $\Gth$, continuous in  $\ovl \Gth\setminus\{0\}$ and vanishing on $\prt\Gth\setminus\{0\}$.

 \smallskip
 
\noindent{\it Step 1: $\Gth:=\BBR^N_+$}. Then $\Gth^{\ell}=\BBR^N$. Letting $c\to\infty$ in $(\ref{uniq9})$ yields to
 \bel{uniq10} 
T_\ell[u_{\infty,0}^{\BBR^N_+}]=u_{\infty,0}^{\BBR^N_+}\qq\forall\ell>0.
 \ee
Therefore $u_{\infty,0}^{\BBR^N_+}$ is self-similar and thus under the separable form $(\ref{Sep-sol})$. By \rth{exis-uniq-SS},
 \bel{uniq11}  
u_{\infty,0}^{\BBR^N_+}(x)=|x|^{\frac{q-2}{q-1}}\gw_s(\frac{x}{|x|}).
 \ee

 \smallskip 
 
 \noindent{\it Step 2: $\Gth:=B$ or $B'^c$}. In accordance with our previous notations, we set
 $B^{\ell}=\frac{1}{\ell}B$ and $B'^{c\,\ell}=\frac{1}{\ell}B'^c$ for any $\ell>0$ and we have,
 \bel{uniq11*} T_\ell[u_{\infty,0}^{B}]=u_{\infty,0}^{B^{\ell}} \text{ and }T_\ell[u_{\infty,0}^{B'^c}]=u_{\infty,0}^{B'^{c\,\ell}}
 \ee
 and
  \bel{uniq12*} u_{\infty,0}^{B^{\ell'}}\leq u_{\infty,0}^{B^{\ell}} \leq u_{\infty,0}^{\BBR^N_+}\leq  u_{\infty,0}^{B'^{c\,\ell}} 
 \leq  u_{\infty,0}^{B'^{c\,\ell''}}\qq\forall \,0<\ell\leq\ell',\ell''\leq 1.
\ee
When $\ell\to 0$ $u_{\infty,0}^{B^{\ell}}\uparrow \underline u_{\infty,0}^{\BBR^N_+}$ and $u_{\infty,0}^{B'^{c\,\ell}}\downarrow \overline u_{\infty,0}^{\BBR^N_+}$ where $ \underline u_{\infty,0}^{\BBR^N_+}$ and $\overline u_{\infty,0}^{\BBR^N_+}$ are positive solutions of $(\ref{B})$ in $\BBR^N_+$ such that 
  \bel{uniq13} 
 u_{\infty,0}^{B^{\ell}} \leq \underline u_{\infty,0}^{\BBR^N_+}\leq u_{\infty,0}^{\BBR^N_+}\leq \overline u_{\infty,0}^{\BBR^N_+}\leq  u_{\infty,0}^{B'^{c\,\ell}}\qq\forall\,0<\ell\leq 1.
\ee
This combined with the monotonicity of  $u_{\infty,0}^{B^{\ell}}$ and $u_{\infty,0}^{B'^{c\,\ell}}$ implies that $ \underline u_{\infty,0}^{\BBR^N_+}$ and $\ovl u_{\infty,0}^{\BBR^N_+}$ vanish on $\prt\BBR^N_+\setminus\{0\}$ and are continuous in 
$\overline{\BBR^N_+}\setminus\{0\}$. Furthermore there also holds for $\ell,\ell'>0$,
  \bel{uniq14} 
T_{\ell'\ell}[u_{\infty,0}^{B}]=T_{\ell'}[T_{\ell}[u_{\infty,0}^{B}]]=u_{\infty,0}^{B^{\ell\ell'}} \text{ and }
T_{\ell'\ell}[u_{\infty,0}^{B'^c}]=T_{\ell'}[T_{\ell}[u_{\infty,0}^{B'^c}]]=u_{\infty,0}^{B'^{c\,\ell\ell'}}.
\ee
Letting $\ell\to 0$ and using $(\ref{uniq11*})$ and the above convergence, we obtain
  \bel{uniq15} 
\underline u_{\infty,0}^{\BBR^N_+}=T_{\ell'}[\underline u_{\infty,0}^{\BBR^N_+}] \text{ and }
\overline u_{\infty,0}^{\BBR^N_+}=T_{\ell'}[\overline u_{\infty,0}^{\BBR^N_+}].
\ee
Again this implies that $\underline u_{\infty,0}^{\BBR^N_+}$ and $\overline u_{\infty,0}^{\BBR^N_+}$ are separable solutions of $(\ref{B})$ in $\BBR^N_+$ vanishing on $\prt \BBR_+^N\sms\{0\}$ and continuous in $\overline{\BBR^N_+}\setminus\{0\}$. Therefore they coincide with $u_{\infty,0}^{\BBR^N_+}$.
  \smallskip 
 
 \noindent{\it Step 3: End of the proof}. From $(\ref{uniq8})$ and $(\ref{uniq11*})$ there holds
   \bel{uniq16} u_{\infty,0}^{B^{\ell}}\leq T_{\ell}[u_{\infty,0}^{\Gw}]
 \leq  u_{\infty,0}^{B'^{c\,\ell}}\qq\forall \,0<\ell\leq 1.
\ee
Since the left-hand side and the right-hand side of $(\ref{uniq16})$ converge to the same function $u_{\infty,0}^{\BBR^N_+}(x)$, we obtain 
   \bel{uniq17} 
\lim_{\ell\to 0}\ell^{\frac{2-q}{q-1}}u_{\infty,0}^{\Gw}(\ell x)=|x|^{\frac{q-2}{q-1}}\gw_s(\frac{x}{|x|})
\ee
and this convergence holds in any compact subset of $\Gw$. If we fix $|x|=1$, we derive $(\ref{uniq7})$.\qeda\medskip

\noindent\Remark It is possible to improve the convergence in $(\ref{uniq7})$ by straightening $\prt\Gw$ near $0$ (and thus to replace $u_{\infty,0}^{\Gw}$ by a function $\tilde u_{\infty,0}^{\Gw}$ defined in $B_\ge\cap \BBR^N_+$) and to obtain a convergence in $C^1(\overline {S^{N-1}_+})$. \medskip

Combining this result with \rth{gen+M} we derive
\bcor{sub} Assume $1<q<q_c$ and $0\in\prt\Gw$. If $u$ is a positive solution of $(\ref{B})$ with boundary trace $tr_{\prt\Gw}(u)=(\CS(u),\gm(u))=(\{0\},0)$ then $u\geq u_{\infty,0}^{\Gw}.$
\es

The next result asserts the existence of a maximal solution with boundary trace $(\{0\},0)$. 
\bprop{max} Assume $1<q<q_c$ and $0\in\prt\Gw$. Then there exists a maximal solution $U:=U_{\infty,0}^{\Gw}$ of $(\ref{B})$ with boundary trace $tr_{\prt\Gw}(U)=(\CS(U),\gm(U))=(\{0\},0)$. Furthermore
 \bel{max1} 
\lim_{\tiny\BA{c}\Gw \ni x\to 0\\
\frac{x}{|x|}=\gs\in S^{N-1}_+
\EA}|x|^{\frac{2-q}{q-1}}U_{\infty,0}^{\Gw}(x)=\gw_s(\gs),
 \ee
 locally uniformly on $S^{N-1}_+$.
\es
\Proof {\it Step 1: Existence}. Since $1<q<q_c<\frac{N}{N-1}$, there exists a radial separable singular solution of $(\ref{B})$ in $\BBR^N\setminus\{0\}$,
 \bel{max2} 
U_S(x)=\Gl_{N,q}|x|^{\frac{q-2}{q-1}}\q\text{with }\;\Gl_{N,q}=\left(\frac{q-1}{2-q}\right)^{q'}\left(\frac{(2-q)(N-(N-1)q)}{(q-1)^2}\right)^{\frac{1}{q-1}}.
 \ee
 By \rlemma{estfunct2} there exists $C_4(q)>0$ such that any positive solution $u$ of $(\ref{B})$ in $\Gw$ which vanishes on 
 $\prt\Gw\setminus\{0\}$ satisfies $u(x)\leq C_4(q)|x|^{\frac{q-2}{q-1}}$ in $\Gw$. Therefore,  
 $U^*(x)=\Gl^*|x|^{\frac{q-2}{q-1}}$ with $\Gl^*:=\Gl^*(N,q)\geq \max\{\Gl_{N,q},C_4(q)\}$ is a supersolution of $(\ref{B})$ in $\BBR^N\setminus\{0\}$ and dominates in $\Gw$ any solution $u$ vanishing on $\prt\Gw\setminus\{0\}$. For $0<\ge<\max\{|z|:z\in\Gw\}$, we denote by $u_\ge$ the solution of
  \bel{max4} 
\left\{\BA {ll} 
-\Gd u_\ge+|\nabla u_\ge|^q=0\qq&\text {in }\Gw\setminus B_\ge\\\phantom{-\Gd +|\nabla u_\ge|^q}
u_\ge=0\qq&\text {on }\prt\Gw\setminus B_\ge\\\phantom{-\Gd +|\nabla u_\ge|^q}
u_\ge=\Gl^*\ge^{\frac{q-2}{q-1}}\qq&\text {on }\Gw\cap \prt B_\ge.
\EA\right.\ee
If  $\ge'<\ge$, $u_{\ge'}|_{_{\prt(\Gw\setminus B_{\ge})}}\leq u_{\ge}|_{_{\prt(\Gw\setminus B_{\ge})}}$, therefore
  \bel{max5} 
  u\leq u_{\ge'}\leq u_\ge\leq U^*(x)\qq\text{in }\;\Gw.
  \ee
Letting $\ge$ to zero, $\{u_\ge\}$ decreases and converges to some $U_{\infty,0}^{\Gw}$ which vanishes on  $\prt\Gw\setminus \{0\}$.   By the the regularity estimates already used in stability results, the convergence occurs in $C^1_{loc}(\overline\Gw\setminus\{0\})$, $U_{\infty,0}^{\Gw}\in C(\overline{\Gw}\setminus\{0\})$ is a positive solution of $(\ref{B})$ and it belongs to $C^2(\Gw)$; furthermore it has boundary trace $(\{0\},0)$ and for any positive solution $u$ satisfying $tr_{\prt \Gw}(u)=(\{0\},0)$ there holds
    \bel{max6} 
 u_{\infty,0}^{\Gw} \leq  u\leq U_{\infty,0}^{\Gw}\leq U^*(x).
  \ee
Therefore $ U_{\infty,0}^{\Gw}$ is the maximal solution.\smallskip

  \noindent{\it Step 2: $\Gw=\BBR^N_+$}. Since 
    \bel{max7} T_\ell[U^*]|_{_{|x|=\ge}}=U^*|_{_{|x|=\ge}}\qq\forall\,\ell>0,
    \ee 
    there holds
      \bel{max8}T_\ell[u_\ge]=u_{\frac{\ge}{\ell}}  \ee
      Letting $\ge\to 0$ yields to $T_\ell[U_{\infty,0}^{\BBR^N_+}]=U_{\infty,0}^{\BBR^N_+}$. Therefore 
      $U_{\infty,0}^{\BBR^N_+}$ is self-similar and coincide with  $u_{\infty,0}^{\BBR^N_+}$.
      \smallskip
  
  \noindent{\it Step 3: $\Gw=B$ or $B'^c$}. We first notice that the maximal solution is an increasing function of the domain. Since $T_\ell[u^\Gth_\ge]=u^{\Gth^\ell}_{\frac{\ge}{\ell}}$ 
where we denote by $u^\Gth_\ge$ the solution of $(\ref{max4})$ in $\Gth\setminus B_\ge$ for any $\ell,\ge>0$ and any domain $\Gth$ (with $0\in \prt\Gth$), we derive as in \rprop{sing}-Step 2, using $(\ref{max8})$ and uniqueness,
 \bel{max9} 
 T_\ell[U_{\infty,0}^{B}]=U_{\infty,0}^{B^{\ell}} \text{ and }\;T_\ell[U_{\infty,0}^{B'^c}]=U_{\infty,0}^{B'^{c\,\ell}}
 \ee
 and
  \bel{uniq12} U_{\infty,0}^{B^{\ell'}}\leq U_{\infty,0}^{B^{\ell}} \leq u_{\infty,0}^{\BBR^N_+}\leq  U_{\infty,0}^{B'^{c\,\ell}} 
 \leq  U_{\infty,0}^{B'^{c\,\ell''}}\qq\forall \,0<\ell\leq\ell',\ell''\leq 1.
\ee
As in \rprop{sing}, $U_{\infty,0}^{B^{\ell}}\uparrow \underline U_{\infty,0}^{\BBR^N_+}\leq U_{\infty,0}^{\BBR^N_+}$ and
$U_{\infty,0}^{B'^{c\,\ell}} \downarrow \overline U_{\infty,0}^{\BBR^N_+}\geq U_{\infty,0}^{\BBR^N_+}$ where  $\underline U_{\infty,0}^{\BBR^N_+}$ and $ \overline U_{\infty,0}^{\BBR^N_+}$ are positive solutions of $(\ref{B})$ in $\BBR^N$ which vanish on $\prt\BBR^N_+\setminus\{0\}$ and endow the same scaling invariance under $T_\ell$. Therefore they coincide with $u_{\infty,0}^{\BBR^N_+}$.      \smallskip
  
  \noindent{\it Step 3: End of the proof}. It is similar to the one of \rprop{sing}.\qeda\medskip
  
  Combining \rprop{sing} and \rprop{max} we can prove the final result
  \bth {UNI} Assume $1<q<q_c$ and $0\in\prt\Gw$. Then $U^\Gw_{\infty,0}=u^\Gw_{\infty,0}$.
  \es
\Proof We follow the method used in \cite[Sec 4]{GV}.\smallskip

\noindent{\it Step 1: Straightening the boundary}. We represent $\prt\Gw$ near $0$ as the graph of a $C^2$ function $\gf$
defined in $\BBR^{N-1}\cap B_R$ and such that $\gf(0)=0$, $\nabla\gf(0)=0$ and 
$$\prt\Gw\cap B_R=\{x=(x',x_N):x'\in \BBR^{N-1}\cap B_R,x_N=\gf(x')\}.$$
We introduce the new variable $y=\Gf(x)$ with $y'=x'$ and $y_N=x_N-\gf(x')$, with corresponding spherical coordinates in $\BBR^N$, $(r,\gs)=(|y|,\frac{y}{|y|})$. If $u$ is a positive solution of $(\ref{B})$ in $\Gw$ vanishing on $\prt\Gw \sms \{0\}$, we set $\tilde u(y)=u(x)$, then a technical computation shows that $\tilde u$ satisfies with ${\bf n}=\frac{y}{|y|}$
\bel{uni1}\BA {l}
r^2\tilde u_{rr}\left(1-2\gf_r\langle{\bf n},{\bf e}_N\rangle+\abs{\nabla\gf}^2\langle{\bf n},{\bf e}_N\rangle^2\right)\\[2mm]
+r\tilde u_{r}\left(N-1-r\langle{\bf n},{\bf e}_N\rangle\Gd\gf-2
\langle\nabla'\langle{\bf n},{\bf e}_N\rangle,\nabla'\gf\rangle
+r\abs{\nabla\gf}^2\langle\nabla'\langle{\bf n},{\bf e}_N\rangle,{\bf e}_N\rangle
\right)\\[2mm]
+\langle\nabla' \tilde u,{\bf e}_N\rangle\left(2\gf_r-\abs{\nabla \gf}^2\langle{\bf n},{\bf e}_N\rangle-r\Gd\gf\right)\\[2mm]
+r\langle\nabla' \tilde u_r,{\bf e}_N\rangle\left(2\langle{\bf n},{\bf e}_N\rangle\abs{\nabla \gf}^2-2\gf_r\right)-2\langle\nabla' \tilde u_r,\nabla'\gf\rangle\langle{\bf n},{\bf e}_N\rangle\\[2mm]
+\abs{\nabla\gf}^2\langle\nabla'\langle\nabla'\tilde u,{\bf e}_N\rangle,{\bf e}_N\rangle
-\frac{2}{r}\langle\nabla'\langle\nabla'\tilde u,{\bf e}_N\rangle,\nabla'\gf\rangle
+\Gd'\tilde u\\[2mm]
+r^2\abs{\tilde u_r{\bf n}+\frac{1}{r}\nabla'\tilde u-(\gf_r{\bf n}+\frac{1}{r}\nabla'\gf)\langle \tilde u_r{\bf n}+\frac{1}{r}\nabla'\tilde u,{\bf e}_N\rangle}^q
=0.
\EA\ee
Using the transformation $t=\ln r$ for $t\leq 0$ and $ \tilde u(r,\gs)=r^{\frac{q-2}{q-1}}v(t,\gs)$, 
we obtain finally that $v$ satisfies
\begin{equation}\label{uni2}\BA {l}
\left(1+\ge_1\right)v_{tt}+\left(N-\frac{2}{q-1}+\ge_2\right)v_{t}
+\left(\gl_{N,q}+\ge_3\right)v+\Gd'v\\[3mm]
\phantom{---}
+\langle\nabla'v,\overrightarrow {\ge_4}\rangle+\langle\nabla'v_t,\overrightarrow {\ge_5}\rangle+
\langle\nabla'\langle \nabla' v,{\bf e}_N\rangle,\overrightarrow {\ge_6}\rangle\\[3mm]\phantom{---}
-\abs{(\frac{q-2}{q-1}v+v_t){\bf n}+\nabla'\tilde v+\langle(\frac{q-2}{q-1}v+v_t){\bf n}+\nabla'\tilde v,{\bf e}_N\rangle\overrightarrow\ge_7}^q=0,
\EA\end{equation}
on $(-\infty,\ln R]\ti S^{N-1}_{+}:=Q_R$ and vanishes on $(-\infty,\ln R]\ti \prt S^{N-1}_{+}$, where $$\gl_{N,q}=\left(\frac{2-q}{q-1}\right)\left(\frac{q}{q-1}-N\right).$$ 
Furthermore the $\ge_j$ are uniformly continuous functions of $t$ and $\gs\in S^{N-1}$ for $j=1,...,7$, $C^1$ for $j=1,5,6,7$ and satisfy the following decay estimates
\begin{equation}\label{uni3}\BA {l}
|\ge_j(t,.)|\leq Ce^t\q\text{for }\,j=1,...,7\,\text{ and }\,|\ge_{j\,t}(t,.)|+|\nabla'\ge_j|\leq c_{17}e^t\q\text{for }\,j=1,5,6,7.
\EA\end{equation}
Since $v$, $v_t$ and $\nabla'v$ are uniformly bounded and by standard regularity methods of elliptic equations \cite[Lemma 4.4]{GV}, there exist a constant $c'_{17}>0$ and $T<\ln R$ such that
\begin{equation}\label{uni4}\BA {l}
\norm {v(t,.)}_{C^{2,\gg}(\overline {S^{N-1}_{+}})}+\norm {v_t(t,.)}_{C^{1,\gg}(\overline {S^{N-1}_{+}})}
+\norm {v_{tt}(t,.)}_{C^{0,\gg}(\overline {S^{N-1}_{+}})}\leq c_{17}'
\EA\end{equation} 
for any $\gamma\in (0,1)$ and $t\leq T-1$. Consequently the set of functions $\{v(t,.)\}_{t\leq 0}$ is relatively compact in the $C^2(\overline {S^{N-1}_{+}})$ topology and there exist $\eta$ and a subsequence $\{t_n\}$ tending to $-\infty$ such that $v(t_n,.)\to \eta$ when $n\to\infty$ in $C^2(\overline {S^{N-1}_{+}})$. \smallskip

\noindent{\it Step 2: End of the proof}. Taking $u=u^\Gw_{\infty,0}$ or   $u=U^\Gw_{\infty,0}$, with corresponding $v$, we already know that $v(t,.)$ converges to $\gw_s$, locally uniformly on $S^{N-1}_{+}$. Thus $ \gw_s$ is the unique element in the limit set of $\{v(t,.)\}_{t\leq 0}$ and $\lim_{t\to-\infty}v(t,.)=\gw_s$ in $C^2(\overline {S^{N-1}_{+}})$. This implies in particular 
\begin{equation}\label{uni5}\BA {l}
\lim_{x\to 0}\frac{u^\Gw_{\infty,0}(x)}{U^\Gw_{\infty,0}(x)}=1
\EA\end{equation} 
and uniqueness follows from the maximum principle.\qeda\medskip

As a consequence we have a full characterization of positive solution with an isolated boundary singularity
\bcor{charact} Assume $1<q<q_c$, $0\in\prt\Gw$ and $u\in C(\overline\Gw\setminus\{0\})\cap C^2(\Gw)$ is a nonnegative solution of $(\ref{B})$ vanishing on $\prt\Gw\setminus\{0\}$. Then either there exists $c\geq 0$ such that $u=u_{c\gd_0}$, or 
$u=u^\Gw_{\infty,0}=\lim_{c\to\infty}u_{c\gd_0}$.
\es

\mysection{The supercritical case}
In this section we consider the case $q_c \leq q <2$.
\subsection{Removable isolated singularities}
\bth{Re-p} Assume $q_c \leq q<2$, $0\in\prt\Gw$ and $u\in C(\overline\Gw\setminus\{0\})\cap C^2(\Gw)$ is a nonnegative solution of $(\ref{B})$ vanishing on $\prt\Gw\setminus\{0\}$. Then $u\equiv 0$.
\es
\noindent\Proof {\it Step 1: Integral estimates}. We consider a sequence of functions $\gz_n\in C^\infty(\BBR^N)$ such that $\gz_n(x)=0$ if $|x|\leq \frac{1}{n}$, $\gz_n(x)=1$ if $|x|\geq \frac{2}{n}$, $0\leq \gz_n\leq 1$ and $|\nabla\gz_n|\leq c_{18}n$, $|\Gd\gz_n|\leq c_{18}n^2$ where $c_{18}$ is independent of $n$. As a test function we take $\gx\gz_n$ (where $\gx$ is the solution to $(\ref{eta})$) and we obtain
\bel{est3}\BA {l}
\myint{\Gw}{}\left(|\nabla u|^q\gx\gz_n -u\gz_n\Gd\gx\right) dx=\myint{\Gw}{}u\left(\gx\Gd \gz_n+2\nabla\gx.\nabla\gz_n \right)dx
=I+II.
\EA\ee
Set $\Gw_n=\Gw\cap \{x:\frac{1}{n}<|x|\leq \frac{2}{n}\}$, then $|\Gw_n|\leq c_{18}'(N)n^{-N}$, thus
$$I\leq c_{18}C_4(q)\myint{\Gw_n}{}n^{\frac{2-q}{q-1}+2}\gx dx\leq c_{18}''n^{\frac{2-q}{q-1}+2-1-N}=c_{18}''n^{\frac{1}{q-1}-\frac{1}{q_c-1}}
$$
since $\gx(x)\leq c_3d(x)$. Notice that $\frac{1}{q-1}-\frac{1}{q_c-1}\leq 0$.
$$II\leq c_{18}C_4(q)\myint{\Gw_n}{}n^{\frac{2-q}{q-1}+1}|\nabla\gx| dx\leq c_{19}n^{\frac{2-q}{q-1}+1-N}=c_{19}n^{\frac{1}{q-1}-\frac{1}{q_c-1}}.
$$
Since the right-hand side of $(\ref{est3})$ remains uniformly bounded, it follows from monotone convergence theorem that
\bel{est4}\BA {l}
\myint{\Gw}{}\left(|\nabla u|^q\gx +u\right) dx<\infty.
\EA\ee
More precisely, if $q>q_c$, $I+II$ goes to $0$ as $n\to\infty$ which implies 
$$\myint{\Gw}{}\left(|\nabla u|^q\gx +u\right) dx=0.
$$
Next we assume $q=q_c$. Since $\abs{\nabla u} \in L^{q_c}_d(\Gw)$, $v:=\BBG^{\Gw}[|\nabla u|^{q_c}] \in L^1(\Gw)$. Furthermore, $u+v$ is positive and harmonic in $\Gw$. Its boundary trace is a Radon measure and since the boundary trace $Tr(v)$ of $v$ is zero, there exists $c\geq 0$ such that $Tr(u)=c\gd_0$. Equivalently, $u$ solves the problem
	\bel{est5} \left\{ \BA {ll}
	-\Gd u+|\nabla u|^{q_c}=0\qquad&\text {in }\Gw\\[2mm]
	\phantom{-\Gd +|\nabla u|^{q_c}}
	u=c\gd_0&\text {in }\prt\Gw.
	\EA \right.\ee
Furthermore, since $u\in L^1(\Gw)$, $u(x)\leq cP(x,.)$ in $\Gw$. Therefore, if $c=0$, so is $u$. Let us assume that $c>0$. \smallskip

\noindent {\it Step 2: The flat case}. Assume $\Gw=B_1^+:=B_1\cap\BBR^N_+$. We use the spherical coordinates $(r,\gs)\in [0,\infty)\ti S^{N-1}$ as above. Put
$$\overline f=\myint{S^{N-1}_+}{}f \tl \vgf_1 dS
$$
then
\bel{X1}\overline u_{rr}+\frac{N-1}{r}\overline u_{r}-\frac{N-1}{r^2}\overline u=\overline{\abs{\nabla u}^{q_c}}
\ee
Set $v(r)=r^{N-1}\overline u(r)$, then
\bel{X2}v_{rr}+\frac{1-N}{r}v_{r}=r^{N-1}\overline{\abs{\nabla u}^{q_c}}.
\ee
and
\bel{X3}v_r(r)=r^{N-1}v_r(1)-r^{N-1}\myint{r}{1}\overline{\abs{\nabla u}^{q_c}}(s)ds.
\ee
Since
\bel{X4}\myint{0}{1}r^{N-1}\myint{r}{1}\overline{\abs{\nabla u}^{q_c}}(s)ds=\frac{1}{N}\myint{0}{1}r^{N}\overline{\abs{\nabla u}^{q_c}}(s)ds<\infty
\ee
it follows that there exists $\lim_{r\to 0}v(r)=\ga \geq 0$. By arguing by contradiction, we deduce that $\ga =0$. Hence 
\bel{X7}
\lim_{r\to 0}r^{N-1}\myint{S^{N-1}_+}{}u(r,\gs)\tl \vgf_1(\gs)dS=0.
\ee
By Harnack inequality \rth{harn}, we obtain
\bel{X9*}
\lim_{x\to 0}|x|^N\myfrac{u(x)}{d(x)}=0. \ee
By standard regularity methods, $(\ref{X9*})$ can be improved in order to take into account that $u$ vanishes on $\prt\BBR^N_+\setminus\{0\}$ and we get
\bel{X10}
\lim_{x\to 0}|x|^{N}\frac{u(x)}{d(x)}=0\Longleftrightarrow\lim_{x\to 0}\frac{u(x)}{P^{\BBR^N_+}(x,0)}=0,
\ee
where $P^{\BBR^N_+}(x,0)$ is the Poisson kernel in $\BBR^N_+$ with singularity at $0$. Since $P^{\BBR^N_+}(.,0)$ is a super solution and $u=o(P^{\BBR^N_+}(.,0))$, the maximum principle implies $u=0$.

\smallskip

\noindent {\it Step 3: The general case}. For $\ell>0$, we set
$$v_\ell(x)=T_\ell[u](x)=\ell^{N-1}u(\ell x).
$$
Then $v_\ell$ satisfies
\bel{est-k}\left\{\BA {ll}
-\Gd v_{\ell}+|\nabla v_{\ell}|^{q_c}=0&\qquad\text{in }\Gw^\ell\\[2mm]
\phantom{-\Gd +|\nabla v_{\ell}|^{q_c}}v_{\ell}=0&\qquad\text{on }\prt\Gw^\ell\setminus\{0\}
\EA\right.\ee
Furthermore, $T_\ell[P^{\Gw}]=P^{\Gw^\ell}$ with $P^{\Gw}:=P^{\Gw^1}$ and
$$u(x)\leq cP^\Gw(x,0)\quad\forall x\in \Gw\Longrightarrow v_{\ell}(x)\leq cP^{\Gw^\ell}(x,0)\quad\forall x\in \Gw^\ell.
$$
By standard a priori estimates \cite{Lib}, for any $R>0$ there exists $M(N,q,R)>0$ such that, if $\Gg_R=B_{2R}\setminus B_R$,
\bel{est-k1}\BA {ll}\sup\left\{\abs{v_{\ell}(x)}+\abs{\nabla v_{\ell}(x)}:x\in \Gamma_R\cap\Gw^\ell\right\}\\[2mm]
\phantom{--------}
+\sup\left\{\myfrac{\abs{\nabla v_{\ell}(x)-\nabla v_{\ell}(y)}}{|x-y|^\gg}:(x,y)\in  \Gamma_R\cap\Gw^\ell\right\}\leq M(N,q,R),
\EA\ee
where $\gg\in (0,1)$ is independent of $\ell\in (0,1]$. Notice that these uniform estimates, up to the boundary, hold because the curvature of $\prt\Gw^\ell$ remains uniformly bounded when $\ell\in (0,1]$. By compactness, there exist a sequence $\{\ell_n\}$ converging to $0$ and function $v\in C^1(\overline{\BBR^N_+}\setminus\{0\})$ such that 
$$\sup\left\{\abs{(v_{\ell_n}-v)(x)}+\abs{\nabla (v_{\ell_n}-v)(x)}:x\in \Gamma_R\cap\Gw^{\ell_n}\right\}\to 0
$$
Furthermore $v$ satisfies
\bel{est-k2}\left\{\BA {ll}
-\Gd v+|\nabla v|^{q_c}=0&\qquad\text{in }\BBR^N_+\\[2mm]
\phantom{-\Gd +|\nabla v|^{q_c}}v=0&\qquad\text{on }\prt\BBR^N_+\setminus\{0\}.
\EA\right.\ee
From step 2, $v=0$ and
$$\sup\left\{\abs{v_{\ell_n}(x)}+\abs{\nabla v_{\ell_n}(x)}:x\in \Gamma_R\cap\Gw^{\ell_n}\right\}\to 0;
$$
therefore
\bel{est-k3}
\lim_{x\to 0}|x|^{N-1}u(x)=0\qq \text{and} \qq \lim_{x\to 0}|x|^N\abs{\nabla u(x)}=0.
\ee
Integrating from $\prt\Gw$, we obtain
\bel{est-k4}
\lim_{x\to 0}\myfrac{|x|^{N}}{d(x)}u(x)=0.
\ee
Equivalently $u(x)=o(P^\Gw(x,0))$ which implies $u=0$ by the maximum principle.
\qeda

\subsection{Removable singularities}

The next statement, valid for a  positive solution of
\bel{lin}
-\Gd u=f\qq\text{in }\Gw
\ee
where $f\in L_d^1$, is easy to prove:

\medskip

\bprop {equiv} Let $q>1$ and $u$ be a positive solution of $(\ref{B})$. The following assertions are equivalent:\smallskip

\noindent (i) $u$ is moderate (\rdef{mod}).\smallskip

\noindent (ii)  $u\in L^1(\Gw)$, $\abs{\nabla u} \in L^q_d(\Gw)$.\smallskip

\noindent (iii) The boundary trace of $u$ is a positive bounded measure $\gm$ on $\prt\Gw$.\smallskip
\es
\medskip

Let $\vgf$ be the first eigenfunction of $-\Gd$ in $W_0^{1,2}(\Gw)$ normalized so that $\sup_{\Gw}\vgf=1$ and $\gl$ be the corresponding eigenvalue. We start with the following simple result.\medskip

\blemma{wei} Let $\Gw$ be a bounded $C^2$ domain. Then for any $q\geq 1$, $0 \leq \ga <1$, $\gg\in [0,\gd^*)$ and $u\in C^1(\Gw)$, there holds
\bel{Y1}\BA{l}\myint{\gg<d(x)<\gd^*}{}(d(x)-\gg)^{-\ga}|u|^qdx\\[4mm]\phantom{-----}\leq C_{12}\left((\gd^*-\gg)^{-\ga}\myint{\Gs}{}|u(\gd^*,\gs)|^q dS + \myint{\gg<d(x)<\gd^*}{}(d(x)-\gg)^{q-\ga}|\nabla u|^q dx\right)
\EA\ee
where $C_{12}=C_{12}(\ga,q,\Gw)$. If $1<q<2$ and $u$ is a solution of $(\ref{B})$, we obtain, replacing $d$ by $\vgf$,
\bel{Y2}\myint{\Gw}{}\vgf^{1-q}|u|^q dx\leq C_{13}\left(1+\myint{\Gw}{}\vgf|\nabla u|^q dx\right)
\ee
where $C_{13}=C_{13}(q,\Gw)$.
\es

\noindent\Proof Without loss of generality, we can assume that $u$ is nonnegative. By the system of flow coordinates introduced in section 2.1, for any $x \in \Gw_{\gd^*}$, we can write $u(x)=u(\gd,\gs)$ where $\gd=d(x)$, $\gs=\gs(x)$ and $x=\gs-\gd{\bf n}_\gs$, thus 
$$u(\gd,\gs)-u(\gd^*,\gs)=-\myint{\gd}{\gd^*}\nabla u (\gs-s{\bf n}_\gs). {\bf n}_\gs ds =-\myint{\gd}{\gd^*}\frac{\prt u}{\prt s} (s,\gs)ds,
$$
from which it follows
$$u(\gd,\gs)\leq u(\gd^*,\gs)-\myint{\gd}{\gd^*}\frac{\prt u}{\prt s} (s,\gs)ds.
$$
Thus, multiplying both sides by $(\gd-\gg)^{-\ga}$ and integrating on $(\gg,\gd^*)$,  
\bel{Y3}\BA {l}\myint{\gg}{\gd^*}(\gd-\gg)^{-\ga}u(\gd,\gs) d\gd \\[4mm]\phantom{\myint{\gg}{\gd^*}}
\leq \myfrac{(\gd^*-\gg)^{1-\ga}}{1-\ga}u(\gd^*,\gs) + \myint{\gg}{\gd^*}(\gd-\gg)^{-\ga}\myint{\gd}{\gd^*}\abs{\nabla u(s,\gs)}ds\,d\gd \\[4mm] \phantom{\myint{\gg}{\gd^*}}
=  \myfrac{(\gd^*-\gg)^{1-\ga}}{1-\ga}u(\gd^*,\gs) + \myfrac{1}{1-\ga}\myint{\gg}{\gd^*}(s-\gg)^{1-\ga}\abs{\nabla u(s,\gs)}ds. 
\EA\ee
Integrating on $\Gs$ and using the fact that the mapping is a $C^1$ diffeomorphism, we get the claim when $q=1$. If $q>1$, we apply $(\ref{Y3})$ to $u^q$ instead of $u$ and obtain
	\bel{Y4} \BA {l}\myint{\gg}{\gd^*}(\gd-\gg)^{-\ga}u^q(\gd,\gs) d\gd \\[4mm]\phantom{}
	\leq \myfrac{(\gd^*-\gg)^{1-\ga}}{1-\ga}u^q(\gd^*,\gs) + \myfrac{q}{1-\ga}\myint{\gg}{\gd^*}(s-\gg)^{1-\ga}u^{q-1}\abs{\nabla u(s,\gs)}ds
	\\[4mm]\phantom{}
	\leq \myfrac{(\gd^*-\gg)^{1-\ga}}{1-\ga}u^q(\gd^*,\gs) + \myfrac{q}{1-\ga}\left(\myint{\gg}{\gd^*}(\gd-\gg)^{-\ga}u^q ds\right)^\frac{1}{q'}\left(\myint{\gg}{\gd^*}(\gd-\gg)^{q-\ga}\abs{\nabla u}^qds\right)^\frac{1}{q}.
	\EA\ee
Since the following implication is true
	$$ (A \geq 0, B \geq 0, M \geq 0, A^q \leq M^q + A^{q-1}B) \Lra (A \leq M + B) $$
we obtain
	\bel{Y5} \BA{ll}
	\left(\myint{\gg}{\gd^*}(\gd-\gg)^{-\ga}u^q(\gd,\gs) d\gd\right)^\frac{1}{q} \\[4mm]
	\phantom{------}
	\leq \left(\myfrac{(\gd^*-\gg)^{1-\ga}}{1-\ga}\right)^\frac{1}{q}u^q(\gd^*,\gs) + \myfrac{q}{1-\ga}\left(\myint{\gg}{\gd^*}(\gd-\gg)^{q-\ga}\abs{\nabla u}^qds\right)^\frac{1}{q}.
	\EA\ee
Inequality $(\ref{Y1})$ follows as in the case $q=1$. We obtain $(\ref{Y2})$ with $\gg=0$, $\ga=q-1$ and using the fact that $c_{21}^{-1}d \leq \vgf \leq c_{21}\,d$ in $\Gw$ with $c_{21}=c_{21}(N)$. \qeda\medskip

\medskip


\bth{removth}Assume $q_c \leq q <2$. Let $K\subset\prt\Gw$ be compact such that $C_{\frac{2-q}{q},q'}(K)=0$. Then any positive moderate solution $u\in C^2(\Gw)\cap C(\overline\Gw\setminus K)$ of  $(\ref{B})$ such that $\abs{\nabla u} \in L^q_{d}(\Gw)$ which vanishes on $\prt\Gw\setminus K$  is identically zero.\es
\Proof Let $\eta\in C^2(\Gs)$ with value $1$ in a neighborhood $U_\eta$ of $K$ and such that $0\leq \eta\leq 1$,  consider $\gz=\vgf(\BBP^\Gw[1-\eta])^{2q'}$. It is easy to check that $\gz$ is an admissible test function since $\zeta(x)+|\nabla\zeta (x)|=O(d^{2q'+1}(x))$ in any neighborhood of $\{x \in \prt \Gw: \eta(x)=1\}$. Then
$$\BA {l}\myint{\Gw}{}|\nabla u|^q\gz dx=\myint{\Gw}{}u\Gd\gz dx=-\myint{\Gw}{}\nabla u.\nabla \gz dx.\EA $$
Next
$$\nabla\gz=(\BBP^\Gw[1-\eta])^{2q'}\nabla\vgf-2q'(\BBP^\Gw[1-\eta])^{2q'-1}\vgf\nabla\BBP^\Gw[\eta],
$$
thus
$$\BA {l}\myint{\Gw}{}|\nabla u|^q\gz dx=-\myint{\Gw}{}(\BBP^\Gw[1-\eta])^{2q'}\nabla \vgf.\nabla u dx+2q'\myint{\Gw}{}(\BBP^\Gw[1-\eta])^{2q'-1}\nabla\BBP^\Gw[\eta].\nabla u \,\vgf dx
\\[4mm]\phantom{\myint{\Gw}{}|\nabla u|^q\gz dx}
= \myint{\Gw}{}u\nabla ((\BBP^\Gw[1-\eta])^{2q'}\nabla\vgf)\, dx+2q'\myint{\Gw}{}(\BBP^\Gw[1-\eta])^{2q'-1}\nabla\BBP^\Gw[\eta].\nabla u \,\vgf dx.
\EA$$
Therefore
\bel{Y6*} \BA{ll}
\myint{\Gw}{}(\gl u + \abs{\nabla u}^q)\zeta dx \\ \phantom{--,}
= -2q'\myint{\Gw}{}(\BBP^\Gw[1-\eta])^{2q'-1}u\nabla \vgf . \nabla \BBP^\Gw[\eta]dx + 2q'\myint{\Gw}{}(\BBP^\Gw[1-\eta])^{2q'-1}\vgf\nabla u . \nabla \BBP^\Gw[\eta]dx. \EA \ee
Since $0 \leq \BBP^\Gw[1-\eta] \leq 1$, $\abs{\nabla \vgf} \leq c_{22}$ in $\Gw$ and by H\"older inequality, 
\bel{Y7}\BA {l}
\abs{\myint{\Gw}{}(\BBP^\Gw[1-\eta])^{2q'-1}u\,\nabla \vgf . \nabla\BBP^\Gw[\eta] \, dx}
\leq c_{22}\left(\myint{\Gw}{}\vgf^{1-q}u^{q} dx\right)^{\frac{1}{q}}\left(\myint{\Gw}{}\vgf|\nabla  \BBP^\Gw[\eta]|^{q'}dx\right)^{\frac{1}{q'}}.
\EA\ee
Using $(\ref{Y2})$ and the fact that $\abs{\nabla u} \in L_d^q(\Gw)$, we obtain 
\bel{Y8}\BA {l}
\abs{\myint{\Gw}{}(\BBP^\Gw[1-\eta])^{2q'-1}u\,\nabla \vgf . \nabla\BBP^\Gw[\eta] \, dx}
\leq c_{23}\left(1 + \norm{\nabla u}_{L_d^q(\Gw)}^q\right)^{\frac{1}{q}}\left(\myint{\Gw}{}d|\nabla  \BBP^\Gw[\eta]|^{q'}dx\right)^{\frac{1}{q'}},
\EA\ee
where $c_{23}=c_{23}(N,q,\Gw)$. Using again H\"older inequality, we can estimate the second term on the right-hand side of $(\ref{Y6*})$ as follows   
\bel{Y8*} \BA {l}\myint{\Gw}{}(\BBP^\Gw[1-\eta])^{2q'-1}\vgf \nabla u . \nabla\BBP^\Gw[\eta] \, dx
\leq \left(\myint{\Gw}{}|\nabla u|^q\vgf dx\right)^{\frac{1}{q}}\left(\myint{\Gw}{}\vgf|\nabla \BBP^\Gw[\eta]|^{q'} dx\right)^{\frac{1}{q'}}\\[4mm]
\phantom{\myint{\Gw}{}(\BBP^\Gw[1-\eta])^{2q'-1}\vgf\nabla\BBP^\Gw[\eta].\nabla u \, dx}
\leq c_{21}\norm{\nabla u}_{L^q_d(\Gw)}^q \left(\myint{\Gw}{}d|\nabla \BBP^\Gw[\eta]|^{q'} dx\right)^{\frac{1}{q'}}.
\EA \ee
Combining $(\ref{Y6*})$, $(\ref{Y8})$ and $(\ref{Y8*})$ we derive
\bel{Y9}\myint{\Gw}{}\left(|\nabla u|^q+\gl u\right)\zeta dx
\leq c'_{23}\left(1 + \norm{\nabla u}_{L_d^q(\Gw)}^q\right)^{\frac{1}{q}}\left(\myint{\Gw}{}d|\nabla  \BBP^\Gw[\eta]|^{q'}dx\right)^{\frac{1}{q'}}.
\ee
By \cite[proposition 7' and Lemma 4']{St},
	\bel{Y9'} \myint{\Gw}{}d|\nabla  \BBP^\Gw[\eta]|^{q'}dx \leq c_{24}\norm{\eta}^{q'}_{W^{1-\frac{2}{q'},q'}(\Gs)}=c_{24}\norm{\eta}^{q'}_{W^{\frac{2-q}{q},q'}(\Gs)}, \ee
which implies
	\bel{Y9*} \myint{\Gw}{}\left(|\nabla u|^q+\gl u\right)\zeta dx
	\leq c_{25}\left(1 + \norm{\nabla u}_{L_d^q(\Gw)}^q\right)^{\frac{1}{q}}\norm{\eta}_{W^{\frac{2-q}{q},q'}(\Gs)} \ee
where $c_{25}=c_{25}(N,q,\Gw)$. Since $C_{\frac{2-q}{q},q'}(K)=0$,  there exists a sequence of functions $\{\eta_n\}$ in $C^2(\Gs)$ such that for any $n$, $0\leq \eta_n \leq 1$, $\eta_n \equiv 1$ on a neighborhood of $K$ and $\norm{\eta_n}_{W^{\frac{2-q}{q},q'}(\Gs)}\to 0$ and $\norm{\eta_n}_{L^1(\Gs)} \to 0$ as $n \to \ity$. By letting $n \to \ity$ in $(\ref{Y9*})$ with $\eta$ replaced by $\eta_n$ and $\zeta$ replaced by $\zeta_n:=\vgf (\BBP[1-\eta_n])^{2q'}$, we deduce that $\myint{\Gw}{}\left(|\nabla u|^q+\gl u\right)\vgf dx=0$ and the conclusion follows. \qeda \medskip
\subsection{Admissible measures}
\bth {adm} Assume $q_c\leq q<2$ and let $u$ be a positive moderate solution of $(\ref{B})$ with boundary data $\gm \in \GTM^+(\prt \Gw)$. Then $\gm(K)=0$ for any Borel subset $K\subset\prt\Gw$ such that $C_{\frac{2-q}{q},q'}(K)=0$.
\es
\Proof Without loss of generality, we can assume that $K$ is compact. We consider test function $\eta$ as in the proof of \rth{removth}, put $\gz=(\BBP^\Gw[\eta])^{2q'}\vgf$ and get
\bel{Z1}\myint{\Gw}{}\left(|\nabla u|^q\gz-u\Gd\gz\right)dx=-\myint{\prt \Gw}{}\myfrac{\prt\gz}{\prt\bf n}d\gm.
\ee
By Hopf lemma and since $\eta \equiv 1$ on $K$, 
$$-\myint{\prt \Gw}{}\myfrac{\prt\gz}{\prt\bf n}d\gm\geq c_{26}\gm (K).
$$
Since 
$$ -\Gd\gz=\gl\gz+4q'(\BBP^\Gw[1-\eta])^{2q'-1}\nabla\vgf.\nabla\BBP^\Gw[\eta]
-2q'(2q'-1)(\BBP^\Gw[1-\eta])^{2q'-2}\vgf|\nabla \BBP^\Gw[\eta]|^2,
$$
we get 
\bel{Z2}c_{26}\gm (K)\leq 
\myint{\Gw}{}\left((|\nabla u|^q+u\gl)\gz+4q'(\BBP^\Gw[\eta])^{2q'-1}u\nabla\vgf.\nabla\BBP^\Gw[\eta]
\right)dx.\ee
Using again the estimates $(\ref{Y8})$ and $(\ref{Y9'})$, we obtain as in \rth{removth} 
\bel{Z2*}\BA {l}
\abs{\myint{\Gw}{}(\BBP^\Gw[1-\eta])^{2q'-1}u\,\nabla\BBP^\Gw[\eta].\nabla \vgf \, dx}
\leq c'_{26}\left(1 + \norm{\nabla u}_{L_d^q(\Gw)}^q\right)^{\frac{1}{q}}\norm{\eta}_{W^{\frac{2-q}{q},q'}(\Gs)}.
\EA\ee
Therefore
\bel{Z3} c_{26}\gm (K)\leq \myint{\Gw}{}(|\nabla u|^q+u\gl)\gz dx + c'_{26}\left(1 + \norm{\nabla u}_{L_d^q(\Gw)}^q\right)^{\frac{1}{q}}\norm{\eta}_{W^{\frac{2-q}{q},q'}(\Gs)}. \ee
As in \rth{removth}, since $C_{\frac{2-q}{q},q'}(K)=0$,  there exists a sequence of functions $\{\eta_n\}$ in $C^2(\Gs)$ such that for any $n$, $0\leq \eta_n \leq 1$, $\eta_n \equiv 1$ on a neighborhood of $K$ and $\norm{\eta_n}_{W^{\frac{2-q}{q},q'}(\Gs)}\to 0$ as $n \to 0$. Thus $\norm{\eta_n}_{L^1(\Gs)}\to 0$ and $\zeta_n:=(\BBP^\Gw[\eta_n])^{2q'}\vgf\to 0$ a.e. in $\Gw$. Letting $n \to \ity$ in $(\ref{Z3})$ with $\eta$ and $\zeta$ replaced by $\eta_n$ and $\zeta_n$ respectively and using the dominated convergence theorem, we deduce that $\gm(K)=0$.\qeda

\mysection{The cases $q=1,2$}
For the sake of completeness we present some results concerning the two extreme cases $q=1$, $q=2$.
\subsection{The case $q=2$}
If $u$ is a solution of $(\ref{B})$ with $q=2$, the standard Hopf-Cole  change of unknown $u=\ln v$ shows that $v$ is a positive harmonic
function in $\Gw$. Therefore the boundary behavior of $u$ is completely described by the theory of positive harmonic functions.
The following result is a consequence of the Fatou and Riesz-Herglotz theorems.
\bth{Fatou} Let $u$ be a bounded from below solution of 
\bel{Fa1}\BA {l}
-\Gd u+|\nabla u|^2=0\qquad\text{in }\Gw.
\EA\ee
1- Then there exists $\gf\in L_+^1(\prt\Gw)$ such that for a.e. $y\in\prt\Gw$,
\bel{Fa2}\BA {l}
\displaystyle\lim_{\tiny{\BA{c}x\to y\\
\text{non-tangent.}
\EA}}\!\!\!\!\!\!\!\!\!u(x)=\ln\gf(y).
\EA\ee
2- There exists a positive Radon measure $\gn$ on $\prt\Gw$ such that
\bel{Fa3}\BA {l}
u(x)=\ln \left(\BBP^\Gw[\gn](x)\right)\qquad\forall x\in\Gw.
\EA\ee
\es

\noindent\Remark Formula $(\ref{Fa3})$ implies that $u$ satisfies
\bel{Fa4}\BA {l}
u(x)\leq (1-N)\ln d(x)+c_{27}\qquad\forall x\in\Gw
\EA\ee
for some $c_{27}$ depending on $u$. This implies in particular that $u\in L^1(\Gw)$.
\medskip

In the next result we describe the boundary trace of $u$.

\bprop{Fatou-tr} Let the assumptions of \rth{Fatou} be satisfied and $\gn$ is the boundary trace of $e^u$. Then $u$ admits a boundary trace $tr_{\prt\Gw}(u)=(\CS(u),\gm(u))$. Furthermore\smallskip

\noindent 1- $z\in\CS(u)$ if and only if for every neighborhood $U$ of $z$, there holds
\bel{Fa5}\BA {l}
\displaystyle\lim_{\gd\to 0}\myint{\Gs_\gd\cap U}{}\ln \left(\BBP^\Gw[\gn](x)\right) dS=\infty.
\EA\ee

\noindent 2- $z\in\CR(u)$ if and only if there exists a neighborhood $U$ of $z$, such that
\bel{Fa6}\BA {l}
\displaystyle\sup_{0<\gd\leq \gd_z}\myint{\Gs_\gd\cap U}{}\ln \left(\BBP^\Gw[\gn](x)\right) dS<\infty,
\EA\ee
for some $\gd_z>0$.
\es
\Proof  This is a direct consequence of the Hopf-Cole transformation and of \rprop{L1} and \rth{gen}.\qeda

\bcor{Fatou-cor} Under the assumptions of  \rth{Fatou}, if $\gn\in L^2(\prt\Gw)$, then $\nabla u\in L_d^2(\Gw)$, thus
$\CS(u)=\emptyset$.
\es
\Proof If $\gn\in L^2(\prt\Gw)$, then $\nabla v\in L_d^2(\Gw)$ (see e.g. \cite{St}). Since $u$ is bounded from below by some constant $c$, $v\geq e^c$ and 
$$\myint{\Gw}{}d\abs{\nabla u}^2dx\leq e^{-2c}\myint{\Gw}{}d\abs{\nabla v}^2dx<\infty.
$$
The conclusion follows from \rprop{tra-reg}.\qeda

\subsection{The case $q=1$}
In this paragraph we consider the equation
\bel{hom}\BA {l}
-\Gd u+|\nabla u|=0\qquad\text{in }\Gw.
\EA\ee
Although there is no linearity, the results are of linear type and the properties of bounded from below solutions of $(\ref{hom})$ similar to the ones of positive harmonic functions. Since the nonlinearity $g(|\nabla u|)=|\nabla u|$ satisfies the subcriticality assumption $(\ref{M})$, for any bounded Borel measure $\gm$ on $\prt\Gw$ there exists a weak solution to the corresponding problem $(\ref{M1})$. The following extension of \rth{uni1} holds
\bprop{unix} For any $z\in\prt\Gw$, there exists a unique weak solution $u=u_{\gd_z}$ to 
\bel{hom1}\left\{\BA {lll}
-\Gd u+|\nabla u| &=0\qquad&\text{in }\Gw\\ \phantom{----,,,}
u&=\gd_z\qquad&\text{on }\prt\Gw.
\EA\right.
\ee
\es
\Proof The proof is in some sense close to the one of \rth{uni1} and starts with a pointwise estimate of the gradient of $u$. This estimate is obtained by a different change of scale different to the one of \rlemma{APE-grad}.  With no loss of generality, we can asume $z=0$. For $\ell\in (0,1]$, we set $w_\ell(x)=\ell^{N-1}u(\ell x)$. Then $w_\ell$ satisfies
\bel{hom2}\left\{\BA {lll}
-\Gd w_\ell+\ell|\nabla w_\ell|&=0\qquad&\text{in }\Gw^\ell:=\frac{1}{\ell}\Gw\\ \phantom{-----,,}
w_\ell&=\gd_z\qquad&\text{on }\prt\Gw^\ell.
\EA\right.
\ee
By the maximum principle
\bel{hom3}
0\leq w_\ell(x)\leq \ell^{N-1}P^{\Gw^\ell}(\ell x,0).
\ee
Again the curvature of $\prt\Gw^\ell$ remains bounded as well as the coefficient of $\abs{\nabla w_\ell}$. Therefore an estimate similar to $(\ref{APE-grad3})$ applies under the following form
\bel{hom3*}\BA{lll} \sup\{\abs{\nabla w_\ell(x)}:x\in \Gw^{\ell} \cap (B_2\sms B_{\frac{1}{2}})\} \\
	\phantom{-------}
	\leq c_{28}'\sup\{w_\ell(x): x\in \Gw^{\ell} \cap (B_3\sms B_{\frac{1}{3}})\} \\
	\phantom{-------}
	\leq c_{28}'\ell^{N-1}\sup\{u(\ell x): x\in \Gw^{\ell} \cap (B_3\sms B_{\frac{1}{3}})\} \\
	\phantom{-------}
	\leq c_{29}
	\EA 
\ee
Choosing $\ell x=y$ with $|x|=1$ we derive 
\bel{hom4}
|\nabla u(y)|\leq c_{29}|y|^{1-N}\qquad\forall y\in\Gw.
\ee
The remaining of the proof is similar to the one of \rth{uni1}, with the use of \rlemma{lim} which holds with $q=1$.\qeda
\medskip

The main result concerning the case $q=1$ is the following
\bth {trac}Assume $u$ is a positive solution of $(\ref{hom})$ in $\Gw$, then there exists a bounded positive Borel measure $\gm$ such that $u$ is a weak solution of the corresponding problem $(\ref{M1})$.
\es
\Proof This is a direct consequence of the proof of \rth{gen+M}. If $\CS(u)\neq\emptyset$ and $z$ in $\CS(u)$ there holds
$$u\geq u_{\ell\gd_z}\qquad\forall \ell>0.
$$
Because of uniqueness and homogeneity,  $u_{\ell\gd_z}=\ell  u_{\gd_z}$. Letting $\ell\to\infty$ yields to a contradiction.\qeda
\medskip

\appendix \section{Appendix: Removabibility in a domain}
\setcounter{equation}{0}
In the section we assume that $\Gw$ is a bounded open domain in $\BBR^N$ with a $C^2$ boundary.
\subsection{General nonlinearity}
This appendix is devoted to the following equation
	\bel{A1**} \left\{\BA{ll}
	- \Gd u + g(\abs{\nabla u})= \gn \qq&\text{in }\,\Gw\\[2mm]
	\phantom{- \Gd  + g(\abs{\nabla u})}
	u=0\qq&\text{on }\,\prt\Gw
	\EA\right. 
	\ee
where $g$ is a continuous nondecreasing function vanishing at $0$ and $\gn$ is a Radon measure. By a solution we mean a function $u\in L^1(\Gw)$ such that $g(|\nabla u|)\in L^1(\Gw)$ satisfying
	\bel{AA2} 
	\myint{\Gw}{}\left(-u\Gd\gz+g(\abs{\nabla u})\gz\right)dx=\myint{\Gw}{}\gz d\gn
	 \ee
	 for all $\gz\in X(\Gw)$.
The integral subcriticality condition on $g$ is the following
	\bel{AA3}\int_1^\infty g(s)s^{-\frac{2N-1}{N-1}}ds<\infty
	\ee

\bth{A-gen} Assume $g\in \CG_0$ satisfies $(\ref{AA3})$. Then for any positive bounded Borel measure $\gn$ in $\Gw$ there exists a maximal solution $\ovl u_\gn$ of $(\ref{A1**})$. Furthermore, if $\{\gn_n\}$ is a sequence of positive bounded measures in $\Gw$ which converges to a bounded measure $\gn$ in the weak sense of measures in $\Gw$ and $\{u_{\gn_n}\}$ is a sequence of of solutions of $(\ref{A1**})$ with $\gn=\gn_n$, then there exists a subsequence $\{\gn_{n_k}\}$ such that $\{u_{\gn_{n_k}}\}$ converges to a solution $u_\gn$ of $(\ref{A1**})$ in $L^1(\Gw)$ and $\{g(|\nabla u_{\gn_{n_k}}|)\}$ converges to $g(\abs{\nabla u_\gn})$ in $L^1(\Gw)$.
\es
\Proof Since the proof follows  the ideas of the one of \rth{meas}, we just indicate the main modifications. 

\noindent (i) Considering a sequence of functions $\gn_n\in C^{\infty}_0(\Gw)$ converging to $\gn$, the approximate solutions are solutions of 
	\bel{w_n} \left\{ \BA{lll}
		- \Gd w +g(|\nabla (w + \BBG^{\Gw}[\gn_n])]) = 0 \qq &\text{ in } \Gw\\
		\phantom{- \Gd  +g(|\nabla (w + \BBG^{\Gw}[\gm_n])]) }
		w = 0  &\text{ on } \prt \Gw.
	\EA \right. \ee

\noindent (ii) The convergence is performed using
\bel{R1}
\norm{\BBG^{\Gw}[\gn]}_{L^1(\Gw)}+\norm{\BBG^{\Gw}[\gn]}_{M^{\frac{N}{N-2}}(\Gw)}+\norm{\nabla\BBG^{\Gw}[\gn]}_{M^{\frac{N}{N-1}}(\Gw)}\leq c_1\norm \gn_{\mathfrak M(\Gw)}
\ee
in \rprop{E1}.

\noindent (iii) For the construction of the maximal solution we consider $u_\gd$ solution of

\bel{A-EM2}\left\{\BA {ll}
-\Gd  u_\gd+g(|\nabla  u_\gd|)=\gn\qq&\text{in }\Gw'_\gd\\[2mm]\phantom{-\Gd +g(|\nabla \bar w_\gd|)}
 u_\gd=\BBG^{\Gw}[\gn]\qq&\text{on }\Gs_\gd.
\EA\right.\ee
Then consequently, $0<\gd<\gd'\Longrightarrow  u_\gd\leq  u_{\gd'}$ in $\Gw'_{\gd'}$ and $u_\gd \downarrow \ovl u_\gn$. Using similar arguments as in the proof of \rth{meas} we deduce that $\ovl u_\gn$ is the maximal solution of $(\ref{A1**})$. \qeda


\subsection{Power nonlinearity}
We consider the following equation
	\bel{A1*} - \Gd u + |\nabla u|^q = \gn \ee
where $1<q<2$. The study on the above equation also leads to a critical value $q^*= \frac{N}{N-1}$. In the subcritical case $1<q<q^*$, if $\gn$ is a bounded Radon measure, then the problem
	$$ \left\{ \BA{lll} - \Gd u + \abs{\nabla u}^q &= \gn \qq &\text{ in } \Gw \\
	\phantom{-----,}
	u &= 0 &\text{ on } \prt \Gw 
	\EA \right. $$   
admits a unique solution $u \in L^1(\Gw)$ such that $\abs{\nabla u}^q \in L^1(\Gw)$ (see \cite{BaPo} for solvability of a much more general class of equation). In the contrary, in the supercritical case, an internal singular set can be removable provided that its Bessel capacity is null. More precisely, 
\bth{Re-s-i} Assume $q^* \leq q <2$ and $K\subset \Gw$  is compact. If $C_{1,q'}(K)=0$ then any positive solution $u\in C^2(\overline \Gw\setminus K)$ of
\bel{int1}
-\Gd u+|\nabla u|^q=0
\ee
in $\Gw\setminus K$ remains bounded and can be extended as a solution of the same equation in $\Gw$.
\es 
\Proof Let $\eta\in C_c^\infty(\Gw)$ such that $0\leq\eta\leq 1$, $\eta=1$ in a neighborhood of $K$. Put $\gz=1-\eta$ and take $\gz^{q'}$ for test function, then
	$$-q'\myint{\Gw}{}\gz^{q'-1}\nabla u.\nabla \eta dx-\myint{\prt \Gw}{}\frac{\prt u}{\prt \bf n} dS+\myint{\Gw}{}\gz^{q'}|\nabla u|^q dx=0.
	$$
Since
	$$\abs{\myint{\Gw}{}\gz^{q'-1}\nabla u.\nabla \eta dx}\leq \left(\myint{\Gw}{}\gz^{q'}|\nabla u|^{q} dx\right)^{\frac{1}{q}}\left(\myint{\Gw}{}|\nabla \eta|^{q'} dx\right)^{\frac{1}{q'}}.
	$$
Therefore
	$$\myint{\Gw}{}\gz^{q'}|\nabla u|^q dx\leq \myint{\prt \Gw}{}\frac{\prt u}{\prt \bf n} dS+q'\left(\myint{\Gw}{}\gz^{q'}|\nabla u|^{q} dx\right)^{\frac{1}{q}}\left(\myint{\Gw}{}|\nabla \eta|^{q'} dx\right)^{\frac{1}{q'}},
	$$
which implies
	\bel{ReSI1}\myint{\Gw}{}\gz^{q'}|\nabla u|^q dx\leq  c_{30}\myint{\prt \Gw}{}\frac{\prt u}{\prt \bf n} dS + c_{31}\myint{\Gw}{}|\nabla \eta|^{q'} dx.
	\ee
where $c_i=c_i(q)$ with $i=30,31$. Since $C_{1,q'}(K)=0$, there exists a sequence $\{\eta_n\} \sbs C_c^{\ity}(\Gw)$ such that $0 \leq \eta_n \leq 1$, $\eta_n=1$ in a neighborhood of $K$ and $\norm{\nabla\eta_n}_{L^{q'}(\Gw)}\to 0$ as $n \to \ity$. Then the inequality $(\ref{ReSI1})$ remains valid with $\eta$ replaced by $\eta_n$ and $\gz$ replaced by $\gz_n=1-\eta_n$. Thus, since $\gz_n \to 1$ a.e. in $\Gw$, we get
$$\myint{\Gw}{}|\nabla u|^q dx\leq  c_{30}\myint{\prt \Gw}{}\frac{\prt u}{\prt \bf n} dS.
$$
Hence, from the hypothesis, we deduce that $\abs{\nabla u} \in L^q(\Gw)$. 

Next let $\eta\in C^{\infty}_0(\Gw)$ and $\eta_n$ as above, then
$$\myint{\Gw}{}(1-\eta_n)\nabla \eta.\nabla u dx-\myint{\Gw}{}\eta\nabla \eta_n.\nabla u dx+\myint{\Gw}{}(1-\eta_n)\eta|\nabla u|^q dx=0.
$$
Since $\abs{\nabla u}\in L^q(\Gw)$, we can let $n\to \infty$ and obtain by monotone and dominated convergence
$$\myint{\Gw}{}\left(\nabla \eta.\nabla u +\eta|\nabla u|^q \right)dx=0.
$$
Regularity results  imply that $u\in C^2(\Gw)$.\qeda
\bth{mea-zero} Assume $q^* \leq q <2$ and $\gn\in\mathfrak M^+(\Gw)$. Let $u\in L^1(\Gw)$ with $\abs{\nabla u}\in L^q(\Gw)$ is a solution of $(\ref{A1*})$ in $\Gw$. Then $\gn (E)=0$ on Borel subsets $E\subset \Gw$ such that $C_{1,q'}(E)=0$.
\es
\Proof Since $\gn$ is outer regular, it is sufficient to prove the result when $E$ is compact. Let $\eta_n$ be a sequence as in the previous theorem, then
\bel{W1} \myint{\Gw}{}(\nabla u.\nabla\eta_n+\eta_n|\nabla u|^q)dx=\myint{\Gw}{}\eta_nd\gn\geq\gn (E).
\ee
But the left-hand side of $(\ref{W1})$ is dominated by
$$\left(\myint{\Gw}{}|\nabla \eta_n|^{q'}dx\right)^{\frac{1}{q'}}\left(\myint{\Gw}{}\eta_n|\nabla u|^qdx\right)^{\frac{1}{q}}+\myint{\Gw}{}\eta_n|\nabla u|^qdx,
$$
which goes to $0$ when $n\to\infty$, both by the definition of the $C_{1,q'}$-capacity and the fact that $\eta_n\to 0$ a.e. as $n \to \ity$ and is bounded by $1$. Thus $\gn (E)=0$.\qeda
\end{document}